\def\eqdef{\buildrel \text def \over =}
\def\rest{\mathord{\restriction}}
\def\ms{\medskip}

\def\phi{\varphi}

\def\sat{\models}
\def\su{\subseteq}
\def\a{\alpha}
\def\b{\beta}
\def\g{\gamma}

\def\d{\delta}
\def\l{\lambda}
\def\k{\kappa}
\def\z{\zeta}
\def\th{\theta}

\def\om{\omega}

\def\lng{\langle}
\def\rng{\rangle}
\def\ov{\overline}
\def\sm{\setminus}
\def\cont{{2^{\aleph_0}}}

\def\nac{{\text {nacc}\, }\,}
\def\nacc{\nac}
\def\acc{{\text {acc}\,}\,}
\def\cl{{\text {cl}\,}}

\def\cf{{\text {cf}\,}}

\def\otp{{\text {otp}\,}\,}
\def\ran{{\text  {ran}\,}}

\def\id{{\text {id}}}

\def\tcf{{\text {tcf}}}
\def\pcf{{\text {pcf}}}

\def\On{{\text {On}}}
\def\Reg{{\text{Reg}}}
\def\eub{{\text{eub}}}
\def\drop{{\text {drop}}}
\def\Drop{{\text {Drop}}}
\def\fill{{\text {fill}}}
\def\Fill{{\text {Fill}}\,}

\def\imply{\Rightarrow}

\def\dcsn{\medskip\noindent{\bf Discussion}: }

\outer\long\def\ignore#1\endignore{}

\newcount\itemno
\def\itm{\advance\itemno1 \item{(\number\itemno)}}

\def\aitm{\advance\itemno1
\item{(\letter\itemno)}}

\def\letter#1{\ifcase#1 \or a\or b\or c\or d\or e\or f\or g\or h\or
i\or j\or k\or l\or m\or n\or o\or p\or q\or r\or s\or t\or u\or v\or
w\or x\or y\or z\else\toomanyconditions\fi}
\def\raitm{\advance\itemno1 \item{(\rletter\itemno)}}
\def\rletter#1{\ifcase#1\or `\or a\or b\or c\or d\or e\or f\or g\or
h\or i\or k\or k\or l\or n\or n\or q\or r\or t\or v
\else\toomanyconditions\fi}

\documentclass[12pt]{article}
\usepackage{amstext}
\usepackage{amssymb}
\usepackage{amsfonts}
\usepackage{amsthm}
\usepackage{amsopn}

\newtheorem{theorem}{Theorem}
\newtheorem{lemma}{Lemma}

 \newtheorem{definition}[lemma]{Definition}
  \newtheorem{corollary}[lemma]{Corollary}
 \newtheorem{question}[lemma]{Question}
 
 \newtheorem{claim}[lemma]{Claim}
 \newtheorem{fact}[lemma]{Fact}
 \newtheorem{exercise}{Exercise}
\newtheorem{remark}[lemma]{Remark}

\def\witb{\ms\noindent{\tt Where is this in the book?}\ms}
\def\notneeded{\ms \noindent{\tt The rest of this Section is not
needed for later sections}\ms}

\def\stackunder#1#2{{#2}_{#1} }

\def\dbigcup{\bigcup}
\def\dbigcap{\bigcap}

\begin{document}

\author{Menachem Kojman
\\ Carnegie-Mellon University}
\title{The A,B,C of pcf: a companion to pcf theory, part I}
\date{November 1995}
\maketitle

\section{Introduction}

This paper is intended to assist the reader learn, or even better, 
teach a course in pcf theory. Pcf theory can be described as the
journey from the function $\beth$ --- the second letter of the Hebrew
alphabet --- to the function $\aleph$, the first letter. For English
speaking readers the fewer the Hebrew letters the better, of course;
but during 1994-95 it seemed that for a group of 6 post doctoral
students in Jerusalem learning pcf theory from Shelah required knowing
all 22 Hebrew letters, for Shelah was lecturing in Hebrew.

This paper offers a less challenging alternative: learning pcf
in a pretty relaxed pace, with all details provided in every proof, and in
English.  

Does pcf theory need introduction? This is Shelah's theory of reduced
products of small sets of regular cardinals. The most well-known
application of the theory is the bound $\aleph_{\om_4}$ on the power
set of a strong limit $\aleph_\om$, but  other applications to  set
theory, combinatorics, abelian groups, partition calculus and general
topology exist already (and more are on the way).

The essence of pcf theory can be described in a few sentences. The
theory exposes a robust skeleton of cardinal arithmetic and gives an
algebraic description of this skeleton.
Shelah's philosophy in this matter is plain: the exponent function 
 is misleading when used to measure the collection of  subsets of a singular
cardinal. 

The right way to measure the size of, say, $[\aleph_\om]^{\aleph_0}$ is
by $\cf\lng  [\aleph_\om]^{\aleph_0},\su\rng$, the cofinality of the partial
ordering of all countable subsets of $\aleph_\om$ ordered by inclusion. 

The usual
$\aleph_\om^{\aleph_0}$ is obtained from that number by multiplying by
$2^{\aleph_0}$. While  $2^{\aleph_0}$ is ``wild" and cannot be
bounded in ZFC, a great discovery of Shelah's is that  $\cf\lng 
[\aleph_\om]^{\aleph_0},\su\rng$ is bounded in ZFC. This way of looking at 
$\aleph_\om^{\aleph_0}$ separates  chaos from structure: puts the chaotic
exponent $\cont$ aside, and clears the way for the study of the structure of
$[\aleph_\om]^{\aleph_0}$.

Shelah approximates $\cf\lng
[\aleph_\om]^{\aleph_0},\su\rng$ by an interval of regular cardinals,
whose first element is $\aleph_{\om+1}$ and whose last element is
$\cf\lng [\aleph_\om]^{\aleph_0},\su\rng$, and so that every regular
cardinal $\l$ in this interval is the {\em true cofinality} of a
{\em reduced product} $\prod B_\l/J_{<\l}$ of a set $B_\l\su
\{\aleph_n:n<\om\}$ modulo an ideal $J_{<\l}$ over $\om$.

When this reduction is done, the study of  of $[\aleph_\om]^{\aleph_0}$ can
be continued algebraically within {\em pcf theory}: the theory of reduced products
of small sets of regular cardinals. 

This approach to cardinal arithmetic can be thought of as ``algebraic
set theory" in analogy to algebraic topology. The information provided
by this view of cardinal arithmetic is enormous, and influences almost
every branch of mathematics in which the notion of cardinality is
important.

For example, a recent application of pcf theory is
taking a construction, due to M.~E.~Rudin, of a certain topological space of
cardinality
$\aleph_\om^{\aleph_0}$ and re-enacting it on one on those approximations. The
construction goes through and the result is  space with the original
properties whose cardinality is computable in ZFC (see \cite{dowker}).

\subsection {The relation of this paper to Shelah's book}

Shelah's Cardinal arithmetic book (henceforth ``the book'') was
published about a year ago and covers a large part of pcf theory and
its applications. The book reflects the state of Shelah's research in
pcf as it was in 1989 --- the year Shelah ``sealed'' the book. The
theory has advance considerably since --- and, roughly, doubled its
volume.  An important advance in the development of the theory was the
proof in ZFC of the existence of stationary sets in $I[\l]$. With this
theorem the development of the basics of pcf is more transparent.  The
approach taken here is the one taken by Shelah in two pcf courses he
taught in Jerusalem in 1991 and 1995, and which I tried to immitate in
a course I gave at Carnegie-Mellon in 1994/95.  The ideal $I[\l]$ is
used every other page. Thus an initial effort is required to develop
the properties of $I[\l]$ which may look disproportional to the
simplifications it generates, but this is worth the effort, because
the proofs are smoother, more transparent and more informative than
earlier proofs in which $I[\l]$ is not used.

Answering an implicit question: yes, I recommend reading this paper
before reading the book. You will benefit more from reading the book
after already knowing that, for example, generators of pcf always
exist.

All the theorems in this paper are Shelah's, unless otherwise stated.
Not all the material in this paper is contained in the book, though.

Pointers to relevant places in the book and in other existing
presentations of pcf theory will be found  under the macro

\witb

Many of the chapters in the book are not about pcf theory itself,
but present applications of the theory. I believe that those parts of
the book will be more accessible to a reader familiar with pcf theory,
as presented here. 

\subsection{Style of writing}

I decided to write {\em all} the details. I hope you are happy with
this decision; but at the times you are bored with reading --- don't
blame me; it's Jim Baumgartner's fault. I tried to get Jim through a
two-day crash-course in pcf during a meeting we both attended last
summer, and after several hours Jim looked into my eyes and said
grimmly: ``This is not easy, Menachem.  You keep saying that this is
easy, but it's difficult stuff''.

Not any more. When you read you will see. 

The difference in the style of writing between this paper and the book
lies in what literary  critics call ``the implied reader'', this
imaginary person, half way between the author and no-one, to  whom
the  writing is addressed. Shelah's implied reader is very clever, sharp, has
phenomenal memory, wants to know always the {\em most general} formulation
of a theorem and, really, {\em already knows} pcf and needs only to be reminded
of what she knows.

My implied reader is slower, wants to learn the important case  first and
 generalizations only later, can learn one thing at a time, and can
lose hours over a $\z$ mistaken for a $\xi$. She is pretty much like me. 

For my reader's sake I made the rule to never use either of the
following adverbs in this paper: ``clearly",``obviously", and
``easily". Rules are, however, clearly made to be easily broken at
obviously suitable circumstances.

\subsection{Additional material} I included digressions into topics
which are not needed for the development of pcf theory. For instance,
additional club-guessing theorems and applications to the saturation
of ideals. Those sections will be marked by the macro 

\notneeded

 and
and may be skipped in a  first reading. Some of the additional material is
sketched in exercises.

\subsection{What is missing in this version}

At least three sections need to be added to this paper. One about the
structure of $\pcf A$; one about reconstruction a characteristic
function of a model from pcf scales, and the applications to cardinal
arithmetic; and one about smooth or transitive generators.

The first two are written and just need proof reading. Before I write
the third I need to digest \cite{Sh506}.

\subsection{Acknowledgments}

This is always a very pleasant section to write. I thank all the
participants in the pcf course I taught at Carnegie-Mellon 1994/95:
Mike Albert, Matt Bishop, Rami Grossberg, Olivier Lessman, Ric
Statman, Boban Velickovic and Roberto Virga. I benefited from all the
comments they have made (``you write too much on the board"; ``you
don't write enough on the board"; ``you drop letters at the ends of
words" and ``your $\xi$ looks like $\z$"). But more than this, I feel
they shared my excitement over some beautiful mathematics.

I thank Ferna Hartman for typing several handwritten notes, that later
evolved into this text, and Roberto Virga for introducing me to
AMSLaTeX. My enthusiasm about several AMS  fonts is visible in several
proofs which make a trully indulgent use of $\subsetneqq$ and similar
symbols.

Special thanks are extended to Amir Leshem, Sami Lefsches and Charlie
Beaton who allowed me to photocopy their notes from Shelah's pcf
course in 1991 (after my own notes were lost). I am also thankful to
Mirna Dzamonja and Andrzej Roslanowski --- the first for the heroic
effort of taking Hebrew notes from Shelah's lectures on pcf, the other
for sending me copies of those notes.

Finally, I am more than grateful to Saharon Shelah for creating such wonderful
mathematics, and for teaching it to me.


\section{ Dropping ordinals and guessing clubs}

\begin{definition} For a set $X$ of ordinals let $\acc X$ denote the set of
accumulation points of $X$: $\acc X=\{\a\in On: \a=\sup X\cap \a\}$.
Let $\nacc X\eqdef X\sm \acc X$. Let $\cl X=X\cup \acc X$.
\end{definition}

We prepare some combinatorics to be used later on. We begin with some
very simple operations on ordinals.

\begin{definition}
\label{drop}

\begin{itemize}

\item Let $c\su On$ be a  set of ordinals, and
let $\a\in On$ be an ordinal.  Define $\drop(\a,c)\eqdef \sup c\cap
\a$ if $c\cap \a\not=\emptyset$. If $c\cap \a=\emptyset$
then $\drop(\a,c)$ is not defined.

\item If $X\su On$ and and $c\su On$ are sets of ordinals define $\Drop
(X,c)\eqdef\{\drop(\a,c):\a\in X\}$.
\end{itemize}
\end{definition}
Let us list a few simple facts about the operations of dropping
ordinals into  sets defined above:

\begin{fact}

\label{dropfact}
 Suppose $c$ and $X$ are sets of ordinals and  $\a\in On$.

\begin{enumerate}
\item $\drop (\a,c)\le \a$ or $\drop (\a,c)$ is undefined. Equality
holds if and only if $\a\in \acc c$.

\item If $\drop(\a,c)$ is defined then $\drop (\a,c)\in \cl c$

\item $\Drop(X,c)\su \cl c$

\item If $\a_1\le\a_2$ are ordinals and $\drop(\a_1,c)$ is defined then
$\drop (\a_2,c)$ is defined too and $\drop(\a_1,c)\le \drop(\a_2,c)$.

\item $\a\mapsto\drop(\a,c)$ is a  order homomorphism on an end
segment of $X$ and consequently $\otp \Drop(X,c)\le \otp X$.

\item If $c_1\su c_2$ are sets of ordinals and $\a\in On$ and
$\drop(\a,c_1)$ is defined then also $\drop (\a,c_2)$ is defined and
$\drop (\a,c_1)\le \drop (\a,c_2)$;

\item If $\lng c_i:i<i(*)\rng$ is such that $c_i\su On$ and
$c_i\su c_j$ for all $j\le i<i(*)$ then for every $\a\in On$ there
exists some $i(\a)<i(*)$ such that $\drop(\a,c_i)$ {\it stabilizes} at
$i(\a)$, that is  either $\drop(\a,c_i)$ is undefined or $\drop(\a,c_i)$ is
constant for all $i(\a)\le i<i(*)$.

\item If $\lng c_i:i<\k\rng$ is sequence of sets of
ordinals weakly decreasing in $\su$, $\k$ a regular cardinal and $X\su On$,
$|X|<\k$ then there exists some
$i(X)<\k$ such that $\Drop (X,c_i)$ {\it stabilizes} at $i(X)$, that
is $\Drop(X,c_{i(X)})=\Drop (X,c_i)$ for all $i(X)\le i<\k$.

\item $\acc\Drop(X,c)=\acc X\cap \acc c$.
\end{enumerate}
\end{fact}

\begin{proof}  The first statement follows directly from the definition of
$\drop(\a,c)$, and every other statement except the last one follows
from the earlier ones and  the wellfoundedness of the ordinals.

To see the last one suppose that $\b\in \acc \Drop(X,c)$ and let
$\lng\drop(\a_i,c):i<i(*)\rng$ be strictly increasing and cofinal in $\b$
with $\a_i\in X$ for all $i<i(*)$. For every $i<i(*)$ necessarily
$\a_i\le \drop(\a_{i+1},c)$, or else $\drop (\a_i,c)\ge
\drop(\a_{i+1},c)$.   Since $\drop (\a_i,c)\in \cl c$ for all
$i$ by 2, we conclude that $\b\in \acc(\acc c)=\acc c$. Since
$\drop(\a_i,c)\le\a_i$ by 1 and $\a_i\le \drop(\a_{i+1},c)<\b$ we
conclude that $\b\in \acc X$.

Conversely, suppose that $\b\in\acc X\cap \acc c$. Find a set $\{
\b_i:i<i(*)\}\su \b\cap c$ unbounded in $\b$.  Since $\b\in \acc X$,
the ordinal $\a_i:=\min \{\b\cap X\sm (\b_i+1)\}$ is defined for all
$i<i(*)$. Therefore $\b_i\le \drop (\a_i,c)<\b$. It follows that
$\b\in \acc \Drop(X,c)$.
\end{proof}

Now let us use those operations to prove the existence of ``club  guessing
sequences".

\begin{theorem} \label{cg1} Suppose that $\kappa^+<\l$, that $\k,\l$
are regular cardinals, and $S\su \{\d<\l:\cf \d=\k\}$ is a given stationary
subset of $\l$. Then there is a sequence $\ov C=\lng c_\d:\d\in S\rng$
satisfying
\begin{enumerate}

\item $c_\d\su\d$ is a closed subset of $\d$
\item For every club $E\su \l$ the set $N(E)\eqdef \{\d\in S:c_\d\su
E\,\wedge\, \sup c_\d=\d\,\wedge\,\otp c_\d=\k\}$
is a stationary subset of $\l$.
\end{enumerate}

\end{theorem}

\begin{remark} If $\lng a_\d:\d\in
S\rng$ satisfies $a_\d\su\d$ and condition 2  in the Theorem then
$\lng \cl c_\d\cap \d:\d\in S\rng$ satisfies both 1 and 2.
Secondly, replacing
$c_\d$ by any cofinal subset does not spoil 2.
\end{remark}

\begin{proof} Let us define by induction on $i<\k^+$ a sequence $\ov C_i=\lng
c^i_\d:\d\in S\rng$ and a club $E_i\su \l$ as follows. Let $E_0=\l$
and let $\ov C_0=\lng c^0_\d:\d\in S\rng$ be chosen such that for
every $\d\in S$ the set $c^0_\d$ is a club of $\d$ of order type $\k$.

Two of the  induction hypotheses are: $c^i_\d=\Drop(c^0_\d,E_i)$ is a closed
subset of $\d$ and $j<i\imply E_i\su E_j$. This holds for $i=0$.

If $i<\k^+$ is limit and $E_j,\ov C_j$ are defined for all $j<i$, let
$E_i:=\cap_{j<i}E_j$. Since $i<\l$, this is indeed a club of $\l$. Let
$c^i_\d:=\Drop(c^0_\d,E_i)$ for all $\d\in S$. By Fact
\ref{dropfact}.9 $c^i_\d$ is a closed subset of $\d$, (and by Fact
\ref{dropfact}.5 its order type is $\le\k$.) 

 Suppose that $i=j+1<\k^+$ and that $E_j,\ov C_j$ are defined. If
there exists some club $E\su \l$ such that for all $\d\in E\cap S$
either $c^j_\d$ is not a club of $\d$ or $c_\d\not\su E$, then choose
such $E$ and let $E_i=E_j\cap E$.

If, however, there is no club $E\su \l$ as required in the definition
of $E_i$, then for every club $E\su \l$ the set $N^j(E):=\{\d\in
S:\d=\sup c^j_\d\,\wedge\, c^i_\d\su E\}$, is stationary.  To see this,
suppose that $N^j(E)\cap E'=\emptyset$ for some club $E\su\l$, namely
$\d\in S\cap E'\imply c^j_\d\not\su E$ for some club $E'\su \l$.  Then $E\cap
E'$ is a club of $\l$ which is as required for the definition of $E_i$.

Therefore, in the case that no $E_{j+1}$ can be found, the proof is
done: $c^j_\d=\Drop (c^0_\d,E_j)$ is a closed subset of $\d$ by Fact
\ref{dropfact}(8) and
its order type is $\le\k$ by Fact 1.5, so if unbounded in $\d$ it is a club
of $\d$ of order type $\k$.  And, by the previous paragraphs, the
prediction clause 2 in the theorem holds for $\lng c^j_\d:\d\in S \rng$.

Assume, then, that $E_{i+1}$ is defined for all $j<\k^+$ and we shall
derive a contradiction.  Let $E=\bigcap _{i<\k^+}E_i$.  Since
$\k^+<\l$, this is a club of $\l$.  Let $\d\in \acc E\cap S$. The
sequence $\lng c^i_\d:i<\k^+\rng$ must stabilize at some $i(\d)<\k^+$
by Fact \ref{dropfact}.8, since $\lng E_i:i<\k^+\rng$ decreasing and
$c^i_\d=\Drop(c^0_\d,E_i)$. Now since $\d\in \acc c^0_\d\cap \acc E_i$
for all $i<\k^+$, by Fact \ref{dropfact}.9 we conclude that $c^i_\d$
is a club of $\d$.
Choose $j\ge i(\d)$. Since $\d\in E_{j+1}$, either $c^j_\d$ is not a club
of $\d$ or $c_\d\not\su E_{j+1}$; but $c^j_\d$ is a club of $\d$;
therefore $c^j_\d\not\su E_{j+1}$. However $c^{i+1}_\d\su \cl
E_{j+1}\cap \l=E_{j+1}$. So $c^i_\d\not=c^{j+1}_\d$, contrary to
stabilization.\end{proof}

 We make the following observation: the proof above gives a
slightly stronger theorem:

\begin{theorem}\label{cg:E} Suppose  $\k^+<\l$ and $\k,\l$ are
regular cardinals, and
$S\su\{\d<\l:\cf\d=\k\}$ is a given stationary set. Then for every
sequence
$\lng c_\d:\d\in S\rng$ satisfying that $c_\d\su\d$ is a club of
$\d$ of order type $\k$ and every club $E'\su \l$ there exists a club
$E\su E'$ such that $\ov C:=\lng \Drop (c_\d,E):\d\in S\rng$ is a club
guessing sequence, namely satisfies the conditions of Theorem \ref{cg1}.
\end{theorem}

What this theorem says is that every sequence of the right form,
supported by a stationary set of $\l$, can be made into a club
guessing sequence by dropping each of its members into a club  $E$
which is contained in a prescribed club $E'\su \l$.

If $\k>\aleph_0$ the operation $\Drop(c^0_\d,E_i)$ can be
replaced by $c^0_i\cap E_i$ in the proof of Theorem \ref{cg1}. In this case
the corresponding version of Theorem \ref{cg:relativized} is:

\begin{theorem}\label{cg:relativized} Suppose  $\k^+<\l$ and $\k,\l$ are
regular cardinals, and
$S\su\{\d<\l:\cf\d=\k\}$ is a given stationary set. Then for every
sequence
$\lng c_\d:\d\in S\rng$ satisfying that $c_\d\su\d$ is a club of
$\d$ of order type $\k$ and every club $E'\su \l$ there exists a club
$E\su E'$ such that $\ov C:=\lng c_\d\cap E:\d\in S\rng$ is a club
guessing sequence, namely satisfies the conditions of Theorem \ref{cg1}.
\end{theorem}

 In other words, any  sequence
of the right form becomes a club guessing sequence when {\sl relativized} to
some club of
$\l$. Theorem \ref{cg:relativized} of  course does not work for $\k=\om$.

\notneeded

We introduce now a second operation on ordinals, which is a generalization
of drop.

\ignore
 difference between Drop and Fill can be described by the following
quotation of a friend of mine, married to her third husband. She summed up
her opinion about divorce as follows: ``I don't believe any more in
{\sl replacing}; only in {\sl adding}". Dropping a set of ordinals into,
say, a club, {\sl replaces} the members of the set by smaller ordinal from
the club; Filling just joins new members to the set from the members of the
club. The filling will be done, however, by dropping some sets --- which we
will fix --- into the club.
\endignore

\begin{definition} Let $\l$ be regular uncountable, $\b\in \l$ and and let
$c\su 
\l$.
\begin{enumerate}

\item Fix a sequence $\lng e_\b:\b\in\l\cap \acc \l\rng$ such that $e_\b$ is
a club of $\b$ and $\otp e_b=\cf\b$.

\item let $\fill(\b,c):=\Drop (e_\b,c)$ if $\b$ is limit and else let
$\fill(\b,c):=\emptyset$.

\end{enumerate}
\end{definition}

The operation $\fill$ depends on a choice of  a sequence $\lng e_\b:\b\in
\acc\l\rng$, but since for our purposes here the actual choice of this
sequence does not matter, we do not incorporate it into the notation.

 We shall use $\fill$ to prove
a club guessing theorem on cofinality
$\k$ of
$\k^+$. The guessing in this case is weaker, but still useful (See [KS2] for
an application). Let
$S^\l_\k$ denote the set of all elements of $\l$ whose cofinality is
$\k$.

\begin{theorem} Let $\k>\aleph_0$ be a regular cardinal and suppose
$S\su
\{\d<\k^+:\cf\d=\k\}$ is a given stationary set. There exists a
sequence $\ov C=\lng c_\d:\d\in S\rng$ such that
\begin{enumerate}

\item $c_\d\su\d$ is a club of $\d$ and $|c_\d|=\k$
\item For every club $E\su\k^+$ the set $N(E)=\{\d\in S:
\d=\sup\{(\nacc c_\d)\cap E\cap S^{\k^+}_\k\}\}$ is stationary.
\end{enumerate}
\end{theorem}

\dcsn Any sequence $\ov C$ satisfying (0) will have the property that
for every club $E\su\l$, for stationarily many $\d\in S$ unboundedly
many points of $c_\d$ enter the club $E$ --- just because
$\k>\aleph_0$; so the crux of the theorem is in making those  points
 be non-accumulation points of $c_\d$.  Furthermore, it is required
in the theorem that those non-accumulation points be of cofinality
$\k$. It is this additional requirement which makes the proof work: I
do not know how to prove the theorem without satisfying this
requirement. What is being used in the proof is the fact that
$S^{\k^+}_\k$ {\it does not reflect}: If $\b<\k^+$ is a limit
ordinal, then there is a club $e_\b$ of $\b$ such that $e_\b\cap
S^{\k^+}_\k=\emptyset$.

\begin{proof}  By removing a non-stationary set from $S$ we may assume that
$S\su \acc\acc \k^+$ --- namely that every ordinal in $S$ is a limit
of limit ordinals.

Let $c^0_\d\su\acc\k^+$ be a club of $\d$ for all $\d\in S$ and set $\ov
C_0=\lng c^0_\d:\d\in S\rng$.

 Let us define by induction on $n\le\omega$ a club $E_n\su \l$ and a
sequence $\ov C_n:=\lng c_\d:\d\in S\rng$. We
 fix $\lng e_\b:\b<\k^+\,\wedge\, \b\in \acc\k^+\rng$ such that
$e_\b\su\b$ is a club of $\b$ and $\otp e_\b=\cf\b$ to be used in
the definition of  $\fill$.

Let $E_0=\acc \k^+$; $\ov C_0$ is already chosen. Suppose $E_n,\ov
C_n$ are defined, and that $c^n_\d$ is a club of $\d$ with
$|c_\d|=\k$. If for every club $E\su \k^+$ the set $\{\d\in S:
\d=\sup\{\nacc c_\d\cap E\cap S^{\k^+}_\k\}\}$ is stationary, then
$\ov C_n$ is as required by the theorem, and we are done. Else, find a
club $E\su \l$ such that $\d\in S\cap E\imply \sup\{(\nacc c_\d)\cap
E\cap S^{\k^+}_\k\}<\d$ and let $E_{n+1}:=E_n\cap  E$.

Let $c^{n+1}_\d=c^n_\d\cup \bigcup\{\fill(\b,E_{n+1})\sm
\drop(\b,c^n_\d):
\b\in \nacc c^n_\d\sm (E_{n+1}\cup S^{\k^+}_\k)\}$

\noindent
{\bf Explanation}: we are trying to correct $c^n_\d$ to overcome the
counterexample $E_{n+1}$. If a non-accumulation point $\b$ is of the
wrong cofinality or does not enter $E_{n+1}$, we fill in the interval of
$c^n_\d$ just below $\b$  elements of $E_{n+1}$ which are obtained  by
dropping $e_\b$ into $E_{n+1}$. If $\b\in\nacc c^n_\d$ is in $E_{n+1}$ and
has the right cofinality, we add nothing in the interval below it, and $\b$
remains a non-accumulation point also in $c^{n+1}_\d$.

Since $c^n_\d$ is a club of $\d$, by Fact \ref{dropfact}.9 also
$c^{n+1}_\d$ is a club of $\d$.

Suppose that $E_n$ and $\ov C_n$ are defined for all $n<\om$ and we
shall derive a contradiction. Let $E'= \bigcap_nE_n$ and let
$E=\acc (E'\cap S^{\k^+}_\k)$. $E\su \k^+$ is a club of
$\k^+$. Fix $\d\in S\cap E$.

Since $\d\in E_{n+1}$ for all $n$, $\gamma_n:=\sup\{(\nacc c_\d)\cap
E_n\}<\d$. Let $\gamma=\sup\{\gamma_n:n<\om\}$. Since $\cf
\d=\k>\aleph_0$, $\gamma<\d$. Pick $\a\in E'\cap \d$ such that
$\a>\gamma$ and $\cf\a=\k$. Such an $\a$ exists because $\d$ is a
limit of  $E'\cap S^{\k^+}_\k$.

Since $\a>\gamma_0$, $\a\notin \nacc c^0_\d$. But also $\a\notin \acc
c^0_\d$, because $\otp c^0_\d=\k$, and therefore $\a\in \acc
c^0_\d\imply \cf\a <\k$. So $\a\notin c^0_\d$ altogether.

Suppose now that $\a\notin c^n_\d$ and let $\b_n:=\min \{c^n_\d\sm
\a\}$. Thus $\b_n>\a$. Let $\a_n=\sup \a\cap c^n_\d$. since $c^n_\d$ is
closed, $\a_n<\a$. We argue now that $\a\notin c^{n+1}_\d$. The only
way $\a$ can join $c^{n+1}_\d$ is by belonging to $\fill(\b_n,E_{n+1})=\Drop
(e_{\b_n}, E_{n+1})$. Well, $\a\notin \nacc
\fill(\b_n,E_{n+1})\cap (\a_n,\b_n)\su \nacc c^{n+1}_\d$ because
$\a>\gamma\ge
\gamma_{n+1}$ and
$\d\in E_{n+2}$. But --- and this is the point where the
non-reflecting property of $S^{\k^+}_\k$ is used --- also $\a\notin
\acc \Drop(e_{\b_n},E_{n+1})$, because $\acc \Drop
(e_{\b_n},E_{n+1})\su \acc e_{\b_n}$ by Fact \ref{dropfact}(8) and $\acc
e_{\b_n}\cap S^{\k^+}_\k=\emptyset$. Thus $\a\notin c^{n+1}_\d$.

We have proved by induction on $n$, then, that $\a\notin c^n_\d$ for
all $n$. Let $\b_n=\min c^n_\d\sm \a$. We obtain the desired
contradiction by showing that the $\b_n$-s are strictly decreasing.
Since $e_{\b_n}$ is unbounded in $\b_n$, we can choose $\gamma\in
(\a,\b_n)$. Certainly $\drop(\gamma,E_{n+1})\ge\a$, and since
$\a\notin A_{\a_n}$ (since $\a\notin c^{n+1}_\d$) the inequality is
sharp. So $c^{n+1}_\d\cap (\a,\b_n)\not=\emptyset$. Thus
$\b_{n+1}<\b_n$.

 We conclude,
then, that some $\ov C_n$ is as required by the Theorem. \end{proof}{}

The final club-guessing theorem we prove is a combination of the
previous two: it combines guessing clubs in the stronger sense of Theorem 1
AND having all non-accumulation points have large cofinality, as in Theorem
2.  A recent
application of this club-guessing principle is the non-saturation of
the non-stationary ideal over all regular cardinals $\l>\aleph_1$ (see
[GS]) which is sketched below in an exercise.

\begin{theorem} \label{cgcof} Suppose that $\k<\th<\l$ are regular
cardinals, and that
$S\su S^\l_\k$ is stationary.  There exists a sequence $\ov
C=\lng c_\d:\d\in S\rng$ such that

\begin{enumerate}
\item $c_\d\su\d$ is a closed subset of $\d$
\item for every club $E\su \l$ the set $\{\d\in
S:c_\d\su E\,\wedge\, \sup c_\d=\d\,\,\wedge\,\, \otp c_\d=\k\,\,\wedge
\,\, \b\in \nacc c_\d\imply \cf\b\ge \theta\}$ is stationary
\end{enumerate}
\end{theorem}

\begin{proof} Let us denote the set $\{\a<\l:\cf\a\ge \theta\}$ by
$S_{\ge\theta}$, and $\{\a<\l:1<\cf\a<\theta\}$ by $S_{< \theta}$.
Let $S\su S^\l_\k$ be given. Using Theorem 1 we fix a sequence $\ov
C_0=\lng c_\d:\d\in S\rng$ that satisfies conditions 1 and 2
there, which are condition 1 here and condition 2 here from which
the clause ``$\,\,\wedge\,\,\b\in
\nacc_\d\imply \cf\b\ge \theta$'' is omitted. It remains  to
`correct' $\ov C_0$ so that condition 2 here will hold in full.

We assume, by replacing each $c_\d$ by $c_\d\cap \acc \l$, if
 necessary, that
$c_\d\su \acc \l$. There is no loss of generality in doing so,
 because for
every club $E\su \l$ the set $\{\d\in S:c_\d\su E\cap\acc\l\}$ 
is stationary.

 Fix a club $e_\b$ of $\b$ with $\otp
e_\b=\cf\b$ for every limit $\b<\l$, to be used as the required parameter
in the operation  $\fill$.

We define below another operation, $\Fill(c_\d,E)$, using a club $E$
as a parameter, which is designed to take care of the cofinality
requirement. Replacing every $c_\d$ by $\Fill(c_\d,E)$ satisfies that
for almost all $\d\in S$ the non-accumulation points of
$\Fill(c_\d,E)$ have cofinality $\ge\th$; but the order type of
$\Fill(c_\d,E)$ may increase, and the guessing of clubs may be
ruined. The proof goes about showing that if $\Fill$ is performed with
a sufficiently `thin' club, then the guessing property is not
spoiled. The cofinality clause follows then from the guessing
property, and the order-type demand is obtained by no more than minor
cosmetics.

 For a given club
$E\su\acc\l$ and
$\d\in S$ let us define $\Fill(c_\d,E)$ in $\om$
approximations as follows:

\begin{itemize}
\item $c^{E,0}_\d=c_\d$

\item $c^{E,n+1}_\d=c^{E,n}_\d\cup \bigcup\{\fill(\b,E)\sm\drop
(\b,c^{E,n}_\d):\b\in\nacc c^{E,n}_\d\,\wedge\,\cf\b<\th\}$

\item $\Fill(c_\d,E):=\bigcup_{n<\om}c^{E,n}_\d$
\end{itemize}

When $E$ is clear  from context, or fixed, we shall write just $c^n_\d$ for
$c^{E,n}_\d$.  We list a few facts about $\Fill(c_\d,E)$ for a fixed club
$E\su \acc \l$ and $c_\d\su\acc \l$:

\begin{fact}
\begin{enumerate}
\item $\fill (\b,E)$ is defined for all $\b\in c^n_\d$, and all $n<\om$.
\item  $c^n_\d\sm c^0_\d\su E$ for all $n<\om$ and $\d\in S$. Therefore 
If $c_\d\su
E$ then $\Fill(c_\d,E)\su E$.

\item $c^n_\d$ is closed for all $\d\in S$ and all $n<\om$
\item If $\b\in \nacc c^n_\d\cap \acc E$ for some $n<\om$, $\d\in S$
 and
$\cf\b<\th$, then $\b\in\acc c^{n+1}_\d$
\end{enumerate}
\end{fact}

The first  Fact is true for $c^0_\d=c_\d$ because $c_\d\su \acc\l$.
For
$n+1$ use the second Fact and the fact that $E\su \acc E$.

The second fact is proved by induction on $n$, using Fact 1.3 and the
definition of $c^{n+1}_\d$. 
The third fact follows by induction from Fact \ref{dropfact}.9.

To prove   Fact 4 suppose $\b\in \nacc c^n_\d\cap \acc E$ and that
$\cf\b<\th$. Forming $c^{n+1}_\d$ from $c^n_\d$ the definition adds 
into the interval $(\drop (\b,c^n_\d),\b)$ all points of
$\Drop(e_\b,E)\sm\drop(\b,c^n_\d)+1$. Since $\b\in \acc E$, $\b\in
\acc\Drop(e_\b,E)$ by Fact \ref{dropfact}.9 and hence $\b\in \acc c^{n+1}_\d$.

Fact 4 is the important one. We are trying to `get rid' of
non-accumulation points whose cofinality is smaller than $\th$. And Fact (4)
above tells us that if the `bad' point $\b\in \nacc c^n_\d\cap S_{<\th}$ 
happens to lie in  $\acc E$, then it becomes an accumulation point  of the
next approximation $c^{n+1}_\d$.

Let us see next that if for some club $E\su \l$ the sequence $\lng
\Fill(c_\d,E):\d\in S\rng$ has  club guessing property, then it also
satisfies the requirement on cofinality of not accumulation points, and
that this is sufficient for the theorem. 

\begin{claim} Suppose there exists a club $E\su\l$  such that  $\{\d\in
S:\Fill(c_\d,E)\su E'\}$ is stationary for every club $E'\su\l$. Then we can
find a sequence
$\ov C=\lng c'_\d:\d\in S\rng$ as required in the Theorem.
\end{claim}

 \begin{proof} [Proof of Claim]

 Suppose a club $E$  is fixed such that $\{\d\in S:\Fill(c_\d,E)\su E'\}$
is stationary for every club $E'\su\l$. Denote the last set by
$N(E')$. We first show that also $N'(E'):=\{\d\in
N(E'):\b\in\nacc c_\d\imply \cf\b\ge\th\}$ is stationary.

 To do
this we make the following simple observation: since $N(E')$ is
stationary for {\sl every} club $E'$, the set   $N(E'\cap \acc
E)$ is stationary. Let $\d\in N(E'\cap \acc E)\su N(E')$. By the definition of
$N(E')$ we know that $\Fill (c_\d,E)\su (E'\cap
\acc E)$. 

We show now that all  points of $\Fill(c_\d,E)$
whose cofinality is smaller than $\th$ are accumulation points. To do so
pick any $\b\in \Fill(c_\d,E)$  with $\cf\b<\th$, and let $n$ be the first
such that $\b\in c^n_\d$. Now $\b\in \acc E$, and if  $\b\in\nacc c^n_\d$
then by Fact (3) above $\b\in\acc c^{n+1}_\d\su\acc c_\d$. Thus, all
non-accumulation points of $\Fill (c_\d,E)\cap S_{<\th}$ are of cofinality
$\ge\th$.

\ms

Finally we thin-out the sets $\Fill(c_\d,E)$ so make them have the right
order type,  without spoiling the fact that their non-accumulation
points have large cofinality: For every $\d\in S$ such that  
$\b\in\Fill(c_\d,E)\imply
\cf\b\ge\th$, find a set
$A_\d\su \nacc \Fill(c_\d,E)$ such that $\otp A_\d=\k$ and $\sup A_\d=\d$
and define
$c'_\d:=\cl A_\d$. Let $c_\d=\emptyset$ for all
other $\d\in S$. Let us show that $\ov C'=\lng c'_\d:\d\in S\rng$ is as
required by the Theorem. Suppose that $E'\su\l$ is a given club. 
We have shown that the set $N'(E')=\{\d\in S:\Fill(c_\d,E)\su
E'\,\wedge\,\b\in \nacc \Fill(c_\d,E)\imply\cf\b\ge\th\}$ is stationary.
For every $\d\in N'(E')$ the set $A_\d$ is contained in $E'$, cofinal in
$\d$, has order-type $\k$ and is contained in $S_{\ge\th}$.   Since
$E'$ is closed, also $c'_\d=\cl A_\d\su E'$.  Finally, if $\b\in\nacc c'_\d$
then necessarily $\b\in A_\d$ and thus $\b\in S_{\ge\th}$, as required.
\end{proof}

\ms
We have yet to show that  a club $E\su\l$ as in the
sufficient condition can be found. The proof uses the eventual
stabilization of $\Fill$, and is not much different than the proof of
Theorem \ref{cg1}.

\begin{claim} There exists a club $E\su \l$ such that for every club
$E'\su\l$ the set $\{\d\in S:\Fill (c_\d,E)\su E'\}$ is stationary.
\end{claim}

\begin{proof}[Proof of Claim]

 Define by induction on $i<\th$ a
decreasing sequence $\lng E_i:i<\th\rng$ of clubs  $E_i\su\l$.

Let $E_0=\acc \l$. If $i<\l$ is limit, let
$E_i:=\bigcap_{j<i}E_j$. Suppose that $i=j+1<\th$ and $E_j$ is
defined. If $\lng \Fill(c_\d,E_j):\d\in S\rng$ satisfies the
conclusion of the claim, then we are done. Else
there is some club $E\su \l$ such that $\d\in E\cap S\imply
\Fill(c_\d,E_j)\not\su E$. Let $E_i:=E_j\cap E$. 

Suppose that $E_i$ is defined for all $i<\th$, and pick
$\d\in\bigcap_{i<\th}E_i$. This is possible because $\th<\l$.

We show that
$\Fill(c_\d,E_i)$ is constant on some end-segment
$[j,\th)$. We write for simplicity $c^{i,n}_\d$ for $c^{E_i,n}_\d$ for
$i<\th$. We will show by induction on
$n$ that
$c^{i,n}_\d$ is constant on
$[j(n),\th)$ for some $j(n)<\th$. This suffices, since we then  set
$j:=\sup\{j(n):n<\om\}$, which, by regularity of $\th$ is
$<\th$  and recall that $\Fill(c_\d,E_i)=\bigcup c^{i,n}_\d$.

For $n=0$ there is little to prove: let $j(0)=0$ and recall that
$c^i_0=c_\d$ for all $i<\th$. Suppose that $c^i_n$ is constant on
$[j(n),\th)$. For every $\b\in c^{j(n),n}_\d\cap S_{<\th}$ there
is some $j<\th$ at which $\fill (\b,E_i)=\Drop(e_\b,E_i)$ stabilizes
by Fact \ref{dropfact}.8, because $|e_\b|=\cf\b<\th$ and $E_i$ is decreasing in
$i$.  By regularity of $\th>|c^0_\d|$ find $j(n)\le j(n+1)<\theta$
such that $\fill (\b,E_i)$ is constant on $[j(n+1),\th)$ for all
$\b\in
\nacc c^{j(n),n}_\d\cap S^\l_{<\th}$. If $i\in [j(n+1),\th)$ then
 $c^{i,n}_\d=c^{j(n+1),n}_\d$ by the induction hypothesis.  The 
definition  of $c^{i,n+1}_\d$ as $c^{i,n}_\d\cup
\bigcup\{\fill(\b,E_i)\sm \drop(\b,c^{i,n}_\d):\b\in\nacc
c^{i,n}_\d\,\wedge\,\cf\b<\th\}$ implies that
$c^{i,n+1}_\d=c^ {j(n+1),n+1}_\d$, because $c^{i,n}_\d=c^{j(n+1),n}_\d$ and
$\fill(\b,E_i) =\fill(\b,E_{i(n+1)})$ for all $\b\in \nacc c^{i,n}$.

Let $j(*)<\th$ be fixed, then, so that $\Fill(c_\d,E_i)$ is constant 
on $[j(*),\th)$.  Since $\d\in E_{j(*)+1}$ it follows that
$\Fill(c_\d,E_i)\not\su E_{j(*)+1}$. But
$\Fill(c_\d,E_{j(*)+1})\su E_{j(*)+1}$ by Fact (1) above because
$c_\d\su E_{j(*)+1}$. So
$\Fill(c_\d,E_{j(*)})\not=\Fill(c_\d,E_{j(*)+1})$, contrary to stabilization.
\end{proof}
\end{proof}

\subsection{\bf Guessing ideals and applications}

We make a few more definitions regarding the club-guessing sequences we met
above. Let $\ov C=\lng c_\d:\d\in S\rng$ be any sequence of sets indexed by
a stationary $S\su\l$ for some regular $\l$. For each of the guessing
requirement presented in Theorems 1-4 there corresponds naturally a {\sl
guessing ideal}. The ideal consists of all subsets of $\l$ which `fail to
guess' some club $E\su \l$. For example, the guessing requirement in
Theorem 1 is that $\{\d\in S:c_\d\su E\}$ is stationary. Therefore the
ideal corresponding to this notion of guessing is
$\{A\su\l:\exists(E\su\l)\,
\hbox{\rm club }\,[\d\in A\cap S\cap E\imply c_\d\not\su E\}$. Saying that
$\lng c_\d:\d\in S\rng$ satisfies the conclusion of Theorem 1 is
equivalent to saying that the guessing ideal is a proper ideal. It
concentrates on
$S$ (that is,
$\l\sm S$ is in it) and extends the non-stationary ideal. It is
$\l$-complete, but not always normal. 

We can relax the guessing requirement, and thus shrink the ideal: if ``to
guess" means ``to be included in $E$ except for an initial segment" then
guessing ideal is $\{A\su \l:\exists E\su \l \,\hbox{club}\,(\d\in
E\cap S\cap A\imply c_\d\not\su^* E\}$ and is sets is normal when the
sequence under discussion satisfies the conclusion of  Theorem 1. 
The symbol $A\su^* B$ means that a proper end-segment of $A$ is contained
in $B$.

Shelah denotes by $\id^b(\ov C)$ the guessing ideal for inclusion, by
$\id^a(\ov D)$ the guessing ideal for inclusion modulo initial segment and
by $\id^p(\ov C)$ the ideal of guessing in the sense of Theorem 3, namely
having unboundedly many non-accumulation points enter the club.

\witb
 
Club guessing is introduced in \S1 of Chapter III, in terms of the guessing
ideals. The existence of club guessing sequences is in \S2. The
operation defined here as ``Drop" is denoted by $gl$ there (for ``glue").
Theorem \ref{cgcof} is not (as far as I know) in the book; it's in a different
book: the non-structure theory book \cite{nonstructure} p.189.

\subsection{ Exercises}

\begin{itemize}

\item Prove that $\id^a(\ov C)$ is a normal ideal when
$\ov C$ is the sequence from Theorem \ref{cg1}. 

\noindent 
Hint: If $E_\a$ witnesses that $A_\a\in\id^a(\ov C)$ the $\bigcap_\a
E_\a$ witnesses that $\bigcup_\a A_\a\in \id^a(\ov C)$

\item Prove that $\Fill(c_\d,E)$ in the proof Theorem \ref{cgcof} is a
closed set of ordinals.

\item (Gitik-Shelah) Prove that if $\k^+<\l$  and
$\k>\aleph_0$ then the non-stationary ideal on $S^\l_\k$ is not
$\l^+$-saturated, namely that there are $\l^+$ many stationary subsets of
$S^\l_\k$ whose pairwise intersections are non-stationary.

\noindent Hint: 
\begin{enumerate}
\item Show that there is {\em no} $\lng c_\d:\d\in S^\l_\k\rng$ such that
$c_\d\su \d$, $|c_\d|=\k$ and $\a\in c_\d\imply\cf\a\ge \k^+\rng$
with the property that for every club $E\su \l$ the set $N(E):=\{\d\in
S^\l_\k:c_\d\su^* E\}$ contains a club intersected with $S^\l_\k$.
\noindent Hint: Define a decreasing sequence of clubs $\lng E_n:n<\om\rng$ so
that $E_{n+1}\su\acc E_n\cap N(E_n)$ and consider the first point in the
intersection of all clubs of cofinality $\k$.

\item If you replace ``contains a club intersected with $S^\l_\k$" by
``stationary", the a sequence as above can be found on every
stationary $S\su S^\l_\k$ by Theorem \ref{cgcof}.

\item Use saturation to show that  for every sequence $\ov C$ as in the previous
paragraph, which is supported by $S$, there is a stationary $S'\su S$
such that the restriction $C\rest S'$ satisfies the condition in 1.
(This is not a contradiction yet, because in 1. you use the fact that
$S=S^\l_\k$).
\noindent Hint: Find a decreasing chain of clubs $\lng E_i:i<\l^+\rng$ and a
chain of stationary sets $S_i:i<\l^+\rng$ such that $S_i=N(E_i)$ and
the difference $S_i\sm S_{i+1}$ is stationary. To define $E_{i+1}$ assume that
$S_i$ is not as required. At limits let $E_i$ be the diagonal intersection of
the previous $E_j$-s and show that $S_i$ is contained modulo the non-stationary
ideal in every $S_j$ for $j<i$. If the process continues $\l^+$ steps then
saturation is violated. 

\item use the previous paragraph to find a maximal antichain in $\l/NS$ of
stationary subsets of $S^\k_\l$, each carrying a sequence $C$
exemplifying 2. By saturation assume the sets in the antichain are
pairwise disjoint. Show that the union of all $\ov C$-s satisfies the
condition forbidden by 1.
\end{enumerate}
\end{itemize}

\section{The ideal $I[\l]$}

In this Section we introduce and develop the basic properties of the main
combinatorial too we shall be using in later sections. This is Theorem
\ref{beauty} below, that asserts the existence of a stationsry $S\su
S^\l_\k$ in $I[\l]$ for all regular $k,\l$ such that $\k^+<\l$.

 Let $\lambda =cf\lambda $ be regular and uncountable cardinal. Let $  
I\left[ \lambda \right] $ be an ideal over $\lambda $ defined as follows:

\begin{definition}
Let $S\subseteq \lambda$. Then $S\in I\left[ \lambda \right] $ if and only
if there exists a sequence $\bar P=\left\langle P_\alpha :\alpha <\lambda
\right\rangle $ and a closed and unbounded set $E\subseteq \lambda $ such
that:

\begin{enumerate}
\item  $P_\alpha \subseteq \mathcal{P}\left( \alpha \right) $ and $  
\left| P_\alpha \right| <\lambda $

\item  If $\delta \in E\cap S$ then $\delta $ is singular and there
exists a set $c\subseteq \delta $ such that $\delta =\sup c,$ $otp\,c<\delta 
$ and $\forall _{r<\delta }\;c\cap \gamma \in \stackunder {{\beta <\delta }}{  
\dbigcup }P_\beta .$
\end{enumerate}
\end{definition}

\begin{fact} $I\left[ \lambda \right] $ is a normal ideal.  
\end{fact} 

\begin{proof}[Proof of Fact] Suppose that $S_\a\in I[\l]$ for $\a<\l$ and that
$E_\a,\ov P_\a=\lng P^\a_\b:\b<\l\rng$ witness $S_\a\in I[\l]$. There
is no loss of generality in assuming that $\lng P^\a_\b:\b<\l\rng$ is
increasing.

Let $E$ be the diagonal intersection of $\{E_\a:\a<\l\}$ and let
$S=\{\a<\l:\exists \b<\a\;[\z\in S_\b]\}$ be the diagonal union of the
$S_\a$. Define $P_\a:=\bigcup_{\b< \a}P^\b_\a$.

Suppose that $\d\in E\cap S$. Then there is some $\b<\a$ such that
$\d\in S_\b$.  Since $\d\in E$ and $E=\{\a:\forall \b<\a\;[\a\in
E_\b]\}$ it holds that $\d\in E_\b$. Thus $\d$ is singular and there
is some club $c\su \d$ with $c\cap \g\in
\bigcup _{\epsilon<\d}P^\b_\epsilon$ for all $\g<\d$. Since by the
definition of $P_\epsilon$ it follows that $\bigcup_{\epsilon}
P^\b_\epsilon\su
\bigcup_{\epsilon<\d}P_\epsilon$ we have that $c\cap \g\in
\bigcup_{\epsilon<\d}P_\epsilon$. This shows that $S\in I[\l]$.
\end{proof}

\begin{fact}\label{ctblinil} For every regular uncountable $\lambda $ the set 
$S_0^\lambda =\left\{ \alpha <\lambda :cf\lambda =\aleph _0\right\}
\in I\left[ \lambda \right] $
\end{fact}

\begin{proof} Let $P_\a=[\a]^{<\aleph_0}$. If $\d\in S^\l_0$ let $c\su
\d$ be any cofinal subset of $\a$ with $\otp c=\om_0$. If $\g<\d$ then
$c\cap\g$ is finite and thus belongs to $P_\g$. 
\end{proof}

\begin{theorem} \label {firstIlambda}
If $\lambda =cf\lambda >\aleph _0$ is a regular uncountable cardinal,
then

$S_{<\lambda }^{\lambda ^{+}}=\left\{ \alpha <\lambda ^{+}:cf\alpha <\lambda
\right\} \in I\left[ \lambda ^{+}\right] .$
\end{theorem}

\begin{proof} Fact \ref{ctblinil} above gives the theorem for $\lambda =\omega
_1,$ so we assume $\lambda >\omega _1.$

For every $\alpha <\lambda ^{+}$ let $\left\langle a_\zeta ^\alpha :\zeta
<\lambda \right\rangle $ be a sequence of closed subsets of $\alpha $
such that:.

\begin{enumerate}
\item  $\forall _\zeta <\lambda \;\;a_\zeta ^\alpha \subseteq \alpha $
is closed and $\left| a_\zeta ^\alpha \right| <\lambda $

\item  $\zeta _1<\zeta _2\Rightarrow a_{\zeta _1}^\alpha \subseteq
a_{\zeta _2}^\alpha $ and for limit $\zeta ,$ $a_\zeta ^\alpha =\stackunder
{{\xi <\zeta }}{\dbigcup }a_\xi ^\alpha $

\item  $\stackunder{\zeta <\lambda }{\dbigcup }a_\zeta ^\alpha =\alpha $
\end{enumerate}

Let $P_\alpha =\left\{ a_\zeta ^\beta \cap \gamma :\beta \leq \alpha
,\;\zeta <\lambda ,\;\gamma \leq \alpha \right\} .$ This definition implies
that $P_\alpha \subseteq \mathcal{P}\left( \alpha \right) $ and $\left|
P_a\right| \leq \lambda <\lambda ^{+}.$ Suppose that $\delta \in S_{<\lambda
}^{\lambda ^{+}}$. Without loss of generality, $\delta >\lambda $.   We need
to show that there is a club $c\subseteq \delta $ of $\delta $ with $c\cap
\gamma \in \stackunder{\beta <\delta }{\dbigcup }P_\beta $ for all $\gamma
<\delta $.   Let $c\subseteq \delta $ be any club of $\delta $ of order type $  
\kappa =cf\,\delta $ and so that $\min c>\lambda $.   For every $\theta \in c$
the set $E_\theta =\left\{ \zeta <\lambda :a_\zeta ^\delta \cap \theta
=a_\zeta ^\theta \right\} $ is a club of $\lambda $ (by the usual back and
forth plus continuity argument). Therefore, $E=\stackunder{\theta \in c}{  
\dbigcap }E_\theta $ is also a club of $\lambda $.  

There must be some index $\zeta _0$ for which $c\subseteq a_{\zeta
_0}^\delta $.   Find a point $\zeta \left( *\right) \in E$ such that $\zeta
\left( *\right) >\zeta _0$.   Now $a_{\zeta (*)}^\delta $ is a closed subset
of $\delta ,$ of cardinality $<\lambda $ and because it contains $c$ as a
subset, it is unbounded in $\delta $.   As $\delta >\lambda ,\;otp\,a_{\zeta
(*)}^\delta <\delta $.   Let $\gamma <\delta $ be arbitrary. To show that $  
a_{\zeta (*)}^\delta \cap \gamma \in \stackunder{\beta <\delta }{\dbigcup }  
P_\beta $ it suffices to show that $a_{\zeta (*)}^\delta \cap \theta \in
P_\theta $ for some $\theta \in c$ which is greater than $\gamma ,$ because $  
P_\theta $ is closed under taking initial segments.

But, fixing such $\theta ,$ we have $a_{\zeta (*)}^\delta \cap \theta
=a_{\zeta (*)}^\theta ,$ as $\zeta (*)\in E\subseteq E_\theta $ and thus
belongs to $P_\theta $.  
\end{proof}

This is the main theorem of this Section:

\begin{theorem}\label{beauty}
If $\k,\l$ are regular cardinals and $\kappa ^{+}<\lambda$ then there
is a stationary set $S\subseteq S^\l_\k $ in $I\left[
\lambda
\right] $.  
\end{theorem}

\begin{proof}
If $\l$ is a successor of regular, then the Theorem follows from the
previous Theorem. The first $\l>\k^+$ for which the previous Theorem
does not apply is $\k^{+\om+1}$. So we may assume that $\k^{++}<\l$. 
We remark that in addition to $\l$ successor of singular, there is
another case which the previous theorem does not cover and this one
does: the case $\l$ regular limit, namely weakly inaccessible.

Fix a
club guessing sequence $\ov C= \lng c_\a:\a\in S^{\k^{++}}_\k\rng$
as in Theorem 1 in the previous section.

\noindent{\bf Description of the proof}.

The proof will involve two elementary chains of models of $H(\chi)$,
for some large enough regular $\chi$. The first chain will be used to
define $\lng P_i:i<\l\rng$; the other, to prove that this choice
works. The first chain is going to be an element of every member in
the second. At some point of the proof, though, we shall need some
set which is definable in the {\sl second} chain to belong to some
member of the {\sl first} chain. We use the prediction, or
``guessing, property of a club guessing sequence to obtain this.

Fix an elementary chain $\ov M:=\lng M_i:i\le\l\rng$ of submodels of
$\lng H(\chi),\in,\rng$ (for a large enough regular $\chi$)
satisfying:

\begin{itemize}
\item $||M_i||<\l$
\item $\lng M_j:j\le i\rng\in M_{i+1}$, $i\in M_i$  and $i\su M_i$
\item $\ov C,\l\in M_0$ and $\k^{++}+1\su M_0$

\end{itemize}

Let us define $P_i:=M_i\cap \mathcal P(i)$. By condition $(1)$ we see that
$|P_i|<\l$.

Let $S\su \l$ be the set

\[
S:=\left\{\a<\l:\cf \a=\k\,\wedge\, \exists c \text {\,club
of\,}\a\big[(\forall\gamma<\a) (c\cap \gamma\in \bigcup_{i<\a}P_i
)\big]\right\}
\]

 The sequence $\lng P_i:i<\l\rng$ witnesses that $S\in I[\l]$
according to the definition of $I[\l]$. All that is
left to be shown is that $S$ is stationary.

\ignore
\noindent{\bf Explanation} If indeed there is a stationary $S\su
S^\l_\k$ in $I[\l]$ them $M_0$ should know about it and about some
witness $\ov P$, by elementarity, and therefore $P_i$ as chosen here
must suffice.
\endignore

We prove now that $S$ is stationary. Fix a club $E\su \l$ and we will
show that $E$ meets $S$. By shrinking $E$ we may assume that $M_i\cap
\l=i$ for all $i\in E$; this follows from the fact that
$\{i<\l:M_i\cap \l=i\}$ is a club of $\l$.

 Define a second elementary chain of models of $H(\chi)$, $\ov N:=\lng
N_\zeta:\zeta\le
\k^{++}\rng$ satisfying:

 \begin{itemize}

\item $\ov M,E,\l,\ov C\in N_0$ and $\k^{++}+1\su N_0$
\item $||N_\zeta||=\k^{++}$
\item $\lng N_\xi:\xi \le \zeta\rng\in N_{\zeta +
1}$ for $\zeta <\k^{++}$
\end{itemize}

 Define $f(\zeta):=\sup N_\zeta \cap \l$. The function $f:(\k^{++}+1)\to  \th$
is increasing and continuous. Denote $\th= f(\k^{++}$. By condition
3 in the choice of $\ov N$, we see that $f\rest \zeta\in N_{\zeta+1}$
for every
$\zeta<\k^{++}$.

 We observe that for every $\zeta\le \k^{++}$ the ordinal $f(\zeta)$
belongs to $E$. This is true by elementarity and the fact that $E\in
N_\zeta$ for all $\zeta$: if $\beta\in N_\zeta\cap \l$ is arbitrary,
then $N_\zeta\sat (\exists \gamma\in \l)\left[ \gamma \in E\,\wedge\,
\gamma>\b\right]$. Therefore there exists $\b<\gamma\in E\cap N$ is
and thus $E$ unbounded below
$f(\zeta)$, implying $f(\zeta)\in E$. Thus $\ran f\su E$.

Turn now to the chain $\ov M$ and work in $M_{\th+1}$. 
Use the fact that $\th\in
M_{\th+1}$ (condition (2) in the choice of $\ov M$) and the
elementarity of $M_{\th+1}$ and choose a function $g\in
M_{\th+1}$ such that  $g:\k^{++}\to \th$ increasing and
continuous  with $\th=\sup\ran g$.

Since both $f$ and $g$ are increasing continuous on $\k^{++}$ with
ranges cofinally contained in $\th$, a standard argument allows us to
fix a club$E\su\k^{++}$ such that $f\rest E=g\rest E$.
 
Use the club guessing property of $\ov C$ to fix $\d<S^{\k^{++}}_\k$
such that $c_\d\su E$, $\d=\sup c_\d$ and $\otp c_\d=\k$.  Define
$c:=f[c_\d]=g[c_\d]$. 

Thus $c\su f(\d)$ is a club of $f(\d)$
of order type $\k$. We already know that $f(\d)\in E$; we will show
that $f(\d)\in S$ by showing that $c\cap
\gamma\in\bigcup_{i<\th_\d}P_i$ for every $\gamma<f(\d)$.

 Let, then, $X=c_\d\cap \z$ for some $\z\in c_\d$ be
an initial segment of $c_\d$ and let $Y:=f[X]$ be the corresponding
initial segment of $c$. As
$c_\d$ and
$\z$ belong to $N_0$, we have $X\in N_0\su N_{\zeta+1}$; and since
$f\rest
\z\in N_{\z+1}$ by condition 3 in the choice of $\ov N$ we conclude
that $Y=f[X]\in N_{\z+1}$. 

If we knew that $Y\in M_i$ for {\sl some} $i<\l$ we would be done:
Suppose that there were
some $i<\l$ such that

\[
H(\chi) \sat\exists i<\l\left[Y\in M_i)\right]
\]

Since $Y,\ov M,\l\in M_{\z+1}$ and  $N_{\z+1}\prec H(\chi)$

\[
N_{\zeta+1}\sat\exists i<\l[Y\in M_i]
\]

by elementarity.  Thus there is some $i<f(\z+1)<f(\d)$ such that
$Y\in P_i$, as required.

But why is there such $i<\l$ at all? The reason is the club guessing
sequence, which enables $\ov M$ to ``predict" the set $Y$. Recall that
$Y=f[X]$.  The parameter $X$ in this definition belongs to
$M_0$. But $f$ may not belong to any  $M_i\in\ov M$. However, $f\rest
c_\d=g\rest c_\d$ and so we can define $Y$ using $g$ instead. The
function $g$ belongs to $M_{\th+1}$, and the proof is now complete,
as $Y=g[X]$ is definable in $M_{\th+1}$ and hence $Y\in
M_{\th+1}$.\end{proof}

In the next Theorem we formulate a more convenient, equivalent definition of
$I[\l]$.

\begin{theorem}  $S\in I\left[ \lambda \right] $ if and only if there
is a club $E\subseteq \lambda $ and a sequence $\left\langle c_\alpha
:\alpha <\lambda
\right\rangle$ satisfying:

\begin{enumerate}
\item[(0)]  $c_\alpha \subseteq \alpha $  and 
$\otp c_\alpha <\alpha$

\item[(1)]  $\beta \in c_\alpha \imply \beta $ is a successor
ordinal and $c_\alpha \cap \beta =c_\beta $

\item[(2)]  If $\delta \in S\cap E$ then $\delta $ is singular,   
$\otp c_\delta =\cf\delta $ and $\delta =\sup c_\delta$.  
\end{enumerate}
\end{theorem}

The condition (1) will be referred to as {\em coherence} of the sequence $\lng
c_\a:\a<\l\rng$.

\begin{proof} Suppose such $E$ and $\left\langle c_\alpha :\alpha
<\lambda \right\rangle$ exist as above. Let $P_\alpha =\left\{
c_\alpha\right\}$. If $\d\in S\cap E$ then $c_\d\su \d$ is unbounded
in $\d$ and $\otp c_\d=\cf \d$. If $\g<\d$ then $c\cap \g=c\cap
\epsilon$ where $\epsilon=\min[c_\d-(\gamma+1)]$ and $c_\d\cap
\epsilon=c_\epsilon\in P_\epsilon$. Thus $S\in I[\l]$

To prove the converse assume that $S\in I\left[ \lambda
\right] $ and fix a club $E\subseteq \lambda $ and a sequence $\left\langle
P_\alpha :\alpha <\lambda \right\rangle $ that satisfy the conditions in
definition 1 for $S$.   We need to produce $E$ and $\left\langle c_\alpha
:\alpha <\lambda \right\rangle $ that satisfy conditions (0) - (2) for $S$.

We make first a few assumptions on $P_\alpha $ which can be made to hold
without increasing $\left| P_\alpha \right| ,$ which is $<\lambda $.   Assume
that for some large enough regular $\chi $ and some $\alpha \in M\prec
H\left( \chi \right) ,\,P_\alpha =M\cap \mathcal{P}^{\left( \alpha \right)
}. $ Therefore if $x\in P_\alpha $ then also $otp\;x\in P_\alpha ;$ for
every limit ordinal $\beta \in P_\alpha $ a club $e_\beta \subseteq \beta $
of $otp\;e_\beta =cf_\beta $ is also in $P_\alpha ,$ and $P_\alpha $ is
closed under set subtraction, union, intersection etc.

Next fix an increasing and continuous sequence $\left\langle \gamma
_i:i<\lambda \right\rangle $ with $\gamma _i<\lambda $ limit ordinal for all 
$i<\lambda $ and $\gamma _{i+1}-\gamma _i\geq \left| \stackunder{\alpha \leq
\gamma _i}{\dbigcup }P_\alpha \right| +\left| \gamma _i\right| +\aleph _0$.  
By thinning $\left\langle \gamma _i:i<\lambda \right\rangle $ we may assume
that $\gamma _i\in E$ and that if $\gamma _i\in S$ then $\gamma _i$ is
singular.

Let $F_i:\stackunder{\alpha \leq \gamma _i}{\dbigcup }P_\alpha \times
\,\gamma _i\rightarrow \left\{ \zeta +1:\gamma _i<\zeta <\gamma
_{i+1}\right\} $ be a 1-1 function. $F$ codes every pair $\left( x,\beta
\right) $ for $x\in \stackunder{\alpha \leq \gamma _i}{\dbigcup }P_\alpha $
and $\beta <\gamma _i$ by a successor ordinal in the interval $\left( \gamma
_i,\gamma _{i+1}\right) $.  

We turn now to defining $c_\alpha $ for $\alpha <\lambda $.

\noindent {\textbf{Case 1}}\textbf{.} $\alpha $ is successor. As
all $\gamma _i$'s are limits, there is a last $i$ such that $\gamma
_i<\alpha ,$ and thus $\alpha <\gamma _{i+1}$.   If $\alpha \notin \,$ran $  
F_i, $ then let $c_\alpha =\emptyset $.  

Else, $\alpha =F_i\left( x,\beta \right) $ for some $x\in \stackunder{\alpha
\leq \gamma _i}{\dbigcup }P_\alpha $ and $\beta <\gamma _i$.   If $\beta \geq
\min x$ let $c_\alpha =\emptyset $ again.

The remaining case is $\beta <\min x$.   If $x\in \stackunder{\alpha \leq
\gamma _{i^{\prime }}}{\dbigcup }\,P_\alpha $ for some $i^{\prime }<i$ let $  
c_\alpha =\emptyset $.   Else:

Let 
\[
c_\alpha =\left\{ 
\begin{array}{rrl}
& \text{(i)} & \zeta \in x\text{ and }otp\text{ }x\cap \zeta \text{ belongs
to }e_\beta \\ 
&  &  \\ 
& \text{(ii)} & j<i\text{ is the least s.t. }x\cap \zeta \in \stackunder{  
\alpha \leq \gamma _j}{\dbigcup }P_\alpha \\ 
F_j\left( x\cap \zeta ,\beta \right) : &  &  \\ 
& \text{(iii)} & \text{for all }\xi \in x\cap \zeta \text{ s.t. }otp\text{ }  
\left( x\cap \xi \right) \in e_\beta \\ 
&  & \text{there is }j^{\prime }<j\text{ s.t. }x\cap \xi \in \stackunder{  
\alpha \leq \gamma _{j^{\prime }}}{\dbigcup }P_\alpha
\end{array}
\right\} . 
\]
{\textbf{Case 2}}. $\alpha $ is the limit: \textit{If possible},
find an unbounded  $c_\alpha \subseteq \alpha$ of  order type $\cf\alpha <\alpha$
with every $\beta \in c_\alpha $ a successor ordinal and $c_\alpha
\cap \beta =c_\beta $.   Otherwise, let $c_\alpha =\emptyset $.  

We check the conditions on $\left\langle c_\alpha :\alpha <\lambda
\right\rangle $.  

\noindent {\textbf{Condition (0)}}: $c_\alpha \subseteq \alpha$ and $otp$
$c_\alpha <\alpha $.   For limit $\alpha $ ondition (0) holds by the choice of
$c_\alpha$.   If $\alpha $ is successor, then every closed subset of $\alpha $
has order type $<\alpha $. 

\noindent {\textbf{Condition (1)}}: Suppose that $\gamma \in
c_\alpha $.   Then $\gamma \in c_\alpha $ is  a
successor ordinal.   Let $\theta \in c_\alpha \cap \gamma $.   By the
definition of $c_\alpha $ there is $\xi \in x$ with $otp$ $x\cap \xi \in
e_\beta $ and $j^{\prime }<i$ least s.t. $x\cap \xi \in \stackunder{\alpha
\leq \gamma _{j^{\prime }}}{\dbigcup }P_\alpha $ with $\theta =F_{j^{\prime
}}\left( x\cap \xi ,\beta \right) $.   (These are conditions (i), (ii) in the
definition of $c_\alpha $). By condition (iii), and as $\xi \in x\cap
\zeta ,$ we have $j^{\prime }<j$ such that $x\cap \xi \in \stackunder{\alpha
\leq \gamma _{j^{\prime }}}{\dbigcup },\;P_\alpha $.  

Let $y:=x\cap \zeta $.   Now $\gamma =F_j\left( y,\beta \right) $ with $\beta
<\min \,y$ and $\xi \in y$ with $otp$ $y\cap \xi =$ $otp$ $x\cap \xi \in
e_\beta $.   Also, $j^{\prime }<j$ is the least s.t. $y\cap \xi =x\cap \xi \in 
\stackunder{\alpha \leq \gamma _{j^{\prime }}}{\dbigcup }P_\alpha $ and $  
F_{j^{\prime }}\left( y\cap \xi ,\beta \right) =\theta $.  

This settles (i), (ii) for $c_\gamma $ with $y$ substituted for $x$.  
Also, (iii) holds. We conclude that $\theta \in c_\gamma $.   The converse
is also true. So we have $c_\alpha \cap \gamma
=c_\gamma $.  \medskip\ 

\noindent {\textbf{Condition (2)}}: We have to say first what $  
E^{\prime }$ is. Let $E^{\prime }=\{\gamma (i):i$ is limit\}. Suppose that $  
\sigma \in S\cap E^{\prime }$.   We know that $\sigma $ is singular (because $  
\left\{ \gamma _i:i<\gamma \right\} \subseteq E$) and that there is a set $  
x\subseteq \delta ,$ $otp$ $x<\sigma $ and $\forall _{\gamma <\sigma
}\;x\cap \gamma \in \stackunder{\alpha <\sigma }{\dbigcup }P_\alpha $.  

We may assume that $x\in P_\sigma $.   Let $\beta =\,otp$ $x$.   We have $\beta
\in P_\sigma $ and by subtracting we assume that $\beta <\min x$.  

Let $\delta =\sup \left\langle \gamma _i:i<i\left( *\right) \right\rangle $
so $\delta =\gamma _{i(*)}$.  

$e_\beta \subseteq \beta $ and $e_\beta $ has order type $cf$ $x=\,cf$ $  
\sigma $.   We know that $e_\beta \in P_\beta $.  

Define $y=\left\{ \zeta \in x:\text{ }otp\text{ }\left( x\cap \zeta \right)
\in e_\beta \right\} $.   Then $otp$ $y=$ $otp$ $e_\beta =$ $cf$ $\delta $ and 
$y$ is club of $\delta $.  

Let $h\left( \zeta \right) =\min \left\{ j:x\cap \zeta \in \stackunder{  
\alpha \leq \gamma _j}{\dbigcup }P_\alpha \right\} $ for $\zeta \in y$.   So $  
h:y\rightarrow \left\{ \gamma _i:i<\lambda \right\} $ and $h$ is
non-decreasing, as $P_{\gamma _i}$ is closed under taking initial segments
for all $i$. In addition, $\forall \zeta \in y\;h(\zeta )<\gamma
_{i(*)}=\delta $.   

Let $z=\left\{ \zeta \in y:\forall \xi \in y\cap \zeta \;h(\xi )<h(\zeta
)\right\} $.  

Let $c=\left\{ F_{h(\zeta )}\left( x\cap \zeta ,\beta \right) :\zeta \in
z\right\} $.  

The set $c$ serves as a candidate to be $c_\delta$.  It is unbounded
and of otp $\cf \sigma $.  All that we need to verify is that $\gamma
\in c\Rightarrow \beta$ is successor and $c\cap \gamma =c_\gamma$.   This
involves checking that for $\alpha =F_{h(\zeta )}\left( x\cap \zeta ,\beta
\right)$ for $\zeta \in z$ then $\zeta $ and $j:=h(\zeta )$, it holds
that $c_\alpha =c\cap \alpha$.  
\end{proof}

\witb

The ideal $I[\l]$ is defined in Definition 2.4(5), p.14. Theorem
\ref{firstIlambda} is mentioned in Remark 2.4A as item (2) in the list of items
that ``will not be used and are included for the reader's amusemant". The
important Theorem \ref{beauty} is not in the book! It was discovered by Shelah
in 1990,  after the book was ``sealed". The application of a stationary set in
$I[\l]$ for obtaining least upper bounds and for finding generators for pcf
(Sections
\ref{lubs} and \ref{gens} below) are in the book, though, where a stationary set
in $I[\l]$ is treated as a  {\em sufficient} condition for both; by now we know
that the existence of stationary subsets in $I[\l]$ for $\l$ succesor of
singular or inaccessible is a Theorem of ZFC. . 

\section{Obtaining least upper bounds} \label{lubs}     

\begin{verse}
Lub, lub, lub\\
All you need is lub\\
tat  tatatata\\
All you need is lub, lub\\
Lub is all you need\\
\end{verse}

In this Section we define and discuss the
relations $<_I$, $\le_i$ and $\lneqq_I$ over ordinal functions from an
infinite set $A$, where $I$ is an ideal over $A$. The important
issues is when a sequence of ordinal functions on $A$ which is increasing
 in $<_I$ has a least, and an exact, upper bound.

Let $A$ be an infinite set, and let $I\su \mathcal P(A)$ be an ideal over $A$.
We denote by $I^*$ the dual filter $\{X\su A:A\sm X\in \}$ and by $I^+$
we denote $\mathcal P(A)\sm I$. We may ocasionaloly refer to sets in $I$
as ``measure zero sets", sets in $I^*$ as ``measure 1 sets" and sets in
$I^+$ as ``positive measure sets" or simply ``positive sets". 

Consider the relations
$<_I$,
$\le_I$ and
$\lneqq_I$ over all functions from $A$ to the ordinals defined as follows:

\begin{itemize}
\item $f=_I g$ if and only if $\{a\in A:f(a)\not=g(a)\}\in I$
\item $f\le_I g$ if and only if $\{a\in A:f(a)> g(a)\}\in I$
\item $f<_I g$ if and only if $\{a\in A:f(a)\ge g(a)\}\in I$
\item $f\lneqq_I g$ if and only if $f\le_I g$ and $\{a\in 
A : f(a)<g(a)\}\in I^+$. 

\end{itemize}

We will also use the notation 
\begin{itemize}

\item $f\le g$ if and only if  $\forall a\in
A[f(a)\le g(a)]$

\item $f<g$ if and only if $\forall a\in
A[f(a)<g(a)]$ 

\item $f\lneqq g$ if and only if $f\le g$ and $f\not= g$
\end{itemize}

Call a sequence $\lng f_\a:\a<\rng$ of
functions
$f_\a:A\to
\On$, for
$\a<\l$, {\it increasing in $<_I$} ($\le_I$, $\lneqq_I$) if and only if
$\a<\b<\l\imply f_\a<_I f_\b$ ($f_\a\le_I f_\b$, $f_\a\lneqq_I f_\b$). 

For the rest of the section fix an infinite set $A$ and an ideal $I\su \mathcal
P(A)$.

\begin{definition}
Suppose $\mathcal F=\{f_\a:\a<\a(*)\}$ is a set of ordinal functions on $A$.
\begin{itemize}
\item A function $g:A\to \On$ is an {\em upper bound} of $\mathcal F$
if only if
$f_\a\le_I g$ for every $\a<\a(*)$.
\item A function $g:A\to \On$ is a {\em least  upper bound}  of
$\mathcal F$ if and only if $g$ is an upper bound of $\mathcal F$ and
$g\le_I g'$ for every upper bound $g'$ of
$\mathcal F$. 
\item A function $g:A\to \On$ is an {\em exact upper bound} of $\mathcal F$ iff
$g$ is a least upper bound of $\mathcal F$ and for every $g'<_I g$ there is some
$\a<\a(*)$ such that
$g'<_I f_\a$.
\end{itemize}
\end{definition} 

We wish to find sufficient conditions fo existence of a lub and an eub
for a sequnce $\ov f$ which is  increasing in $<_I$.  The condition we shall
provide is in Theorem \ref{lubexists} below. To state this condition we
need the following definition of ``obedience", which makes sense for
every $\ov C$, but is useful when $\ov C$ has the coherence property.

\begin{definition}
Let $\bar C=\left\langle
c_\alpha :\alpha <\lambda \right\rangle$ be so that $c_\alpha \subseteq
\alpha$.   A sequence $\bar f=\left\langle f_\alpha :\alpha <\lambda
\right\rangle $ of ordinal functions on $A$ {\em obeys} $\bar C$ if
and only if $\forall \alpha <\lambda \;\forall _\beta \in c_\alpha\;
[f_\beta  <f_\alpha ]$.    
\end{definition}

The next lemmat says that if a sequence $\ov C$ to which $\ov f$
obeys has coherence, then certain subsequences of $\ov f$ are increasing in $<$.

\begin{lemma} \label{window}
Suppose $\ov f=\lng f_\a:\a<\l\rng$ obeys $\ov C=\lng c_\a:\a<\l\rng$.
If $\ov C$ satisfies that $\b\in  c_\a\imply c_\b=c_\a\cap \b$
then for every $\a<\l$ the sequence $\lng f_\b:\b\in
 c_\a\rng$ is increasing in $<$. 
\end{lemma}

\begin{proof}
Let $\b,\gamma\in  c_\a$ for some $\a<\l$ and supose $\b<\gamma$. Since
$c_\a\cap \gamma=c_\gamma$ by coherence 
we have that $\b\in c_\gamma$. By obedience of $\ov f$ we have
$f_\b<f_\gamma$.
\end{proof}

 The following Theorem provides a sufficient condition for the
existence of an eub for an increasing sequence in
$<_I$. 

\begin{theorem} \label{lubexists}
Suppose $\left| A\right| \le\kappa$ and $\kappa
^{+}<\lambda =\cf \lambda $. Suppose $S\subseteq S_{\kappa
^{+}}^\lambda =\left\langle \delta <\lambda :\cf \delta =\kappa
^{+}\right\rangle $ is stationary, $S\in I[\l]$ and that $\ov C$
witnesses $S\in I[\l]$, namely

\begin{itemize}
\item $\ov C=\lng c_\a:\a<\l\wedge\cf\a\le\k^+\rng$

\item  $c_\alpha \subseteq \alpha $  and $\otp c_\alpha <\alpha$

\item  $\beta \in c_\alpha \Rightarrow \beta $ is successor and $  
c_\alpha \cap \beta =c_\beta $

\item  $\otp c_\delta =\kappa ^{+}$ and $  
\sup c_\delta =\delta $ for all $\delta \in S$
\end{itemize}
If a sequence of
ordinal functions $\ov f=\left\langle f_\alpha :\alpha
<\lambda \right\rangle $ is increasing in $<_I$  and {\em obeys} $\ov C$,
then $\ov f$ has an exact upper bound.  
\end{theorem}

\begin{proof} 
First we prove the existence of a lub; then we shall show that every  lub
of $\ov f$ is an eub. 

We find a lub by approximations. We start with an upper bound and
decrease it whenever it is not a lub.  At limit stages the obedience
of $\ov f$ to $\ov C$ is used to produce a ``smaller" upper bound from
the the set of values of all previous upper bounds.  The existence of
lub will be obtained once we show that this process must terminate.

By induction on
$\zeta <\kappa ^{+}$ we define a sequence of upper bounds $g_\zeta $
to $\ov f$ so that:

\begin{equation}\label{decr}
  \xi <\zeta <\kappa ^{+}\Rightarrow g_\zeta \lneqq _Ig_\xi  
\end{equation}

Let $g_0\left( a\right) =\sup\left\{ f_\alpha \left( a\right)
+1:\alpha <\lambda \right\}$.  For every $\alpha <\lambda$,
$g_0>f_\alpha$ so $g_0$ is an upper bound of $\ov f$.

 For successor $\zeta +1$ just choose, if possible, $g_{\zeta +1}$
satisfying
(\ref{decr})  above. If this is not possible, $g_\xi $ is a lub, as
required, and the induction terminates.

Suppose that $\zeta<\k^+$ is limit. We shall show that if $g_\xi$ is
defined for all $\xi<\z$, then also $g_\xi$ can be defined (and
therefore the induction does not terminate at a limit stage $<\k^+$). 

We need to define sets $S_\zeta
\left( a\right)$, auxiliary functions $h_\alpha ^\zeta$ and an index $\alpha
_\zeta <\l$ before defining $g_\zeta $.

\begin{itemize}
\item  For $a\in A$ let $S_\zeta \left( a\right) :=\left\{ g_\xi \left(
a\right) :\xi <\zeta \right\} $ .

\item  for $\a<\l$ let $h_\alpha ^\zeta(a) :=  \min \left\{ S_\zeta
\left( a\right) 
\sm f_\alpha \left( a\right) \right\}$.   
\end{itemize}

Here are the facts we need about the functions $h^\zeta_\a$:

\begin{fact}\label{factsweneed}
\begin{enumerate}
\item $h_\alpha ^\zeta \in 
\stackunder{a\in A}{\prod }S_\zeta \left( a\right) $ and $h_\alpha ^\zeta \geq
f_\alpha $ 

\item  $f_\a(a)\le f_\b(a)\imply h^\z_\a(a)\le h^\z_\b(a)$ for $\a,\b<\l$
 
\item for a fixed $\zeta<\k^+$ the sequence $\lng h^\zeta_\a:\a<\l\rng$ is
increasing in $\le_I$ and for every subsequence $\lng f_\a:\a\in c\rng$ of $\ov
f$, 
$c\su \l$, which is increasing in $<$ the subsequence $\lng h^\z_\a:\a\in c\rng$
is increasing in $\le$. 

\item if $\a,\b<\l$ and $h^\z_\a(a) <h^\z_\b(a)$, then $h^\z_\a(a)<f^\z_\b(a)$

\item  if $\xi <\zeta  <\kappa ^{+}$ and  $\alpha <\l$ then  $h^\xi_\a\ge
h^\z_\a$

\end{enumerate}
\end{fact}

The first item follows directly from the definition of $h^\zeta_\a$.
For the second item fix $\a,\b<\l$. If $f_\z(a)\le f_\b(a)$ then
$h^\z_\a(a)=\min\{S_\z(a)\sm f_\a(a)\}\le
\min\{S_\z(a)\sm f_\b(a)\} =h^\z_\b$. The third item
followd directly from the second. The fourth is immediate from the definitions.
Finally, if
$\xi<\z<\k^+$ then
$S_\xi(a)\su S_\z(a)$ for all $a\in A$ and therefore $h^\z_\a(a)=\min\{S_\z(a)\sm
f_\a(a)\}\le \min\{S_\xi(a)\sm f_\a(a)\}=h^\xi_\a$.

The  existence of the index
$\alpha _\zeta <\l$  is provided  by the following: 

\begin{claim} \label{cutstab} There is some $\alpha _\zeta
<\lambda $ such that 
\begin{equation} \label{cutstabeq}
\alpha _\zeta \leq \alpha <\lambda \Rightarrow h_\alpha ^\zeta =_Ih_{\alpha
_\zeta }^\zeta . 
\end{equation}
\end{claim}

\begin{proof}[Proof of Claim]

 The sequence $\lng h^\z_\a:\a<\l\rng$ is increasing in $\le_I$ by 3
above. If the sequence $\ov h:=\lng h^\z_\a:\a<\l\rng$ does not
stabilize modulo $I$, then $\forall
\alpha <\lambda
\exists
\gamma < \l \left[h_\alpha ^\zeta \lneqq h_\gamma
^\zeta
\right]$, and as $\ov h$ is increasing in $\le_I$   
 there is a club
$E\subseteq
\lambda $ such that $ \b\in E\imply  \forall \alpha <\beta
\left[h_\alpha ^\zeta \lneqq_I h_\beta ^\zeta \right] $.  

Use the fact that $S$ is stationary to find $\delta \in S\cap \acc E$.

Now find a subset $c\su  c_\d$  such that
between any two members of $c$ there is at
least one point of $E$.

By Lemma \ref{window} the sequence $\lng f_\a:\a\in c\rng$ is
increasing in $<$, and therefore, by Fact \ref{factsweneed}.3 the
sequence $\lng h^\z_\a:\a\in c\rng$ is increasing in $\le$. Because
between any two members of $c$ there are points of $E$, the sequence
$\lng h^\z_\a:\a\in c\rng$ is increasing in $\lneqq_I$. Since
$f\le g\wedge f\lneqq_I g \imply f\lneqq g$ for all $f,g:A\to \On$, we
conclude that $\lng h^\z_\a:\a\in c\rng$ is increasing in $\lneqq$.
This contradicts the fact that $|S_\a(a)|\le \k$ for all $a\in A$: the
sequence $\lng h^\z_\a(a):\a\in c\rng$ is increasing in $\le$, and
since $|S_\z(a)|\le \k$ stabilizes at $\a_a<\k^+$. Putting
$\a=\sup\{\a_a:\a\in A\}<\k^+$ we ontain $h^\z_\a=h^\z_{\a+1}$,
contrary to $h^\z_\a\lneqq_I f^\z_{\a+1}$.
\end{proof}

Using Claim \ref{cutstab}, fix, at a limit stage $\z<\k^+$ the ordinal
$\a_\z<\l$ and define $g_\z:=h^\z_{\z_\z}$. For every $\a_\z<\a<\l$ we have,
of course, $g_\z=h^\z_\a$, so we may increase $\a_\z$ at will.

\begin{claim}
$g_\z$ is an upper bound of $\ov f$ and $g_\z\lneqq_I g_\xi$ for every
$\xi<\zeta$.
\end{claim}

\begin{proof}[Proof of Claim]  Let $\a<\l$ be given. By increasing $\a$ we
may assume that $\a\ge\a_\z$. Then $f_\a\le h^\z_\a=_I h^\z_{\a_\z}$.
Thus $g_\z$ is a lub of $\ov f$.

 If $\xi<\z$ then $f_{\a_\z}<_I g_\xi$,
because
$g_\xi$ is an upper bound and $\ov f$ increasing in $<_I$. This means that for
measure 1 set of $a\in A$ we have $g_\xi(a)>f_{\a_\z}(a)$. Since $g_\xi(a)\in
S_\z(a)$, for each such $a$ we have $h^\z_{\a_\z}(a)\le g_\xi(a)$.
This establishes that $g_\z\le_I g_\xi$ for all $\xi<\z$. Since $\z$
is limit and $g_\xi \gneqq g_{\xi+1}$ for $\xi<\z$, it follows the
$g_\z\lneqq g_\xi$ for all $\xi<\z$.
\end{proof}

We prove the existence of lub by showing that the induction
must stop before $\kappa ^{+}$.   We just showed the inductive definition
of $g_\z$ does  not break at limit $\zeta <\kappa ^{+}$.   Therefore it must
break at some $\zeta +1<\kappa ^{+}$; and this means $g_\zeta $ is a lub.

\begin{claim} \label{nodeckapplus}
For some $\z<\k^+$ the upper bound $g_{\z+1}$ is not defined. 
\end{claim}

\begin{proof}[Proof of Claim]

Suppose to the contrary that the induction goes through all $\zeta
<\kappa ^{+}$.

For every limit $\zeta <\kappa ^{+}$  we have 
$g_\zeta =h_{\alpha _\zeta }^\zeta =_I h_\alpha ^\zeta \text{ for all }\alpha
\geq \alpha _\zeta$. 
By regularity of $\lambda >\kappa ^{+}$ we find $\alpha \left( *\right)
<\lambda$ such that $g_\zeta =_Ih_{\alpha \left( *\right) }^\zeta$ for all
limit $\zeta <\kappa ^{+}$.  

The sequence $\left\langle h_{\alpha \left( *\right) }^\zeta :\zeta \in
\acc \k^{+}\right\rangle$ stabilizes at some $\zeta \left( *\right) \in
\acc \k^{+}$ by \ref{factsweneed} 5  above. Let $\zeta >\zeta (*)$ be
limit. So  $g_{\zeta ( *)}\gneqq g_\zeta =g_{\zeta (*) }$ ---
contradiciton to stabilization.
\end{proof}

We show next that the lub we found is an eub.

\begin{claim} \label{eub} Let $g$ be a lub of $\ov f$.
If $g'<_Ig$ then $g'<_If_\alpha $ for some $\alpha
<\lambda$.  

\end{claim}

\begin{proof}[Proof of Claim] Let $B_\alpha =\{ a\in A:f_\a(a)> g'(a)\} $.   For
$\alpha <\beta <\lambda $ we have $B_\alpha \subseteq _I B_\beta $,
because $f_\a<_I f_\b$.  If $B_\alpha /I$ does not stabilize, then
there exists a club $E\subseteq \lambda $ so that $\alpha <\beta \in
E\Rightarrow B_\alpha \subsetneqq_I B_\beta$.  Fix $\delta \in S\cap
\acc E$.  Find cofinal $c\subseteq  c_\d$ with members of $E$
between any two members of $c$.   

Since between any two members of $c$ there is a point of $E$, the
sequence $\lng B_\a:\a\in c\rng$ is increasing in $\subsetneqq_I$. But by
Lemma \ref{window} the sequence $\lng f_\a:\a\in c\rng$ is increasing
in $<$, which implies that $\lng B_\a:\a \in c\rng$ is increasing in
$\su$. Since to increase in both $\su$ and $\subsetneqq_I$ means to
increase in $\subsetneqq$, we cnclude that $\lng B_\a:\a\in c\rng$ is
a strictly increasing sequence of subsets of $A$ of length $\k^+$,
which contradicts $|A|\le \k$. 

We assume, then, that $B_\alpha/I $ stabilizes at  $\alpha
\left( *\right) <\lambda $.   If $B_{\alpha \left( *\right) }\in I^*$, then
$f_{\alpha \left( *\right) }>_Ig'$, and we are done. Else,
for all $\alpha \geq \alpha \left( *\right) $  
\[
f_\alpha \rest A\sm B_{\alpha \left( *\right) }\leq _I g'\rest A\sm
B_{\alpha \left( *\right) } 
\]
Define, then, 
\[
g''=g^{\prime }\rest A\sm B_{\alpha \left( *\right) }\cup
g\rest  B_{\alpha \left( *\right) }. 
\]
$g''$ is an upper bound of $\ov f$, and, $g''  
\lneqq g$ since $g'<_Ig$ and $A\sm B _{\alpha \left(
*\right) }\in I^{+}$.   This contradicts  $g$ being a {\em least} upper
bound.

\end{proof}
\end{proof}

Let us discuss the role of $I[\l]$ in this proof. The obedience of
$\ov f$ to $\ov C$ and the coherence of $\ov C$ implies, by Lemma
\ref{window}, that $\ov f$ is increasing in $<$ ``locally'', or on
subsequences of $\ov f$ indexed by $c_\d\in \ov C$. We know that if
$\ov h^\z$ does not stabilize modulo $I$, it increases in $\lneqq_I$
on a club of $\l$; we also know that it's impossible for any
subsequence of $\ov h^\z$ of length $\k^+$ to increase in $\lneqq$,
becuse the range of each $h^\z_\z\in\ov h^\z$ is in contained in $\prod
S_\z(a)$. 

Now since the set $S\in I[\l]$ is stationary, we can ``trap'' an
accumulation point $\d$ of $E$ in $S$ and by interweaving members of
$c_\d$ with members of $E$ find a sequence of length $\k^+$ which
increases in both relations. In other words, because the subseuences
on which $\ov h$ increases in $<$ are indexed by $c$-s that converge
to all ordinal in a {\em stationary} set, they intertwine with any
potential club on which $\ov h^\a$ may increase.

The same argument is repeated in Claim \ref{eub}: if the sequence of
$B_\a$ does not stabilize modulo $I$ it increases on a club. An
accumulation point of the club is trapped in $S$; and a subsequence of
length $\k^+$ is thus found wich increases in both $\su$ and
$\subsetneqq_I$.

The demand of obedience can be waived in Theorem \ref{lubexists} if we
assume that $2^{|A|}<\l$ (see exercise below). The pigeon-hole
principle can be used to prove the Lemma \ref{cutstab} and Claim
\ref{nodeckapplus}.

A stationary set in $I[\l]$ liberates us, then, from assumptions about
the size of $2^{|A|}$, provided the sequence at hand obeys a suitable
$\ov C$. 

\witb

Compare with Lemma 2.3 in Burke-Magidor. See also Lemma 2.1 in Jech. 
The book mentions $I[\l]$ and a stationary set in $I[\l]$ as a
sufficient condition for existence of eub (\S2 of chapter I), but not as 
a theorem
of ZFC. 

\begin{exercise}

\begin{itemize} 
\item Prove that if $|A|=\k\ge\aleph_0$, $I$ an ideal over $A$  and
$\lng f_\a:\a<\l\rng$ is a sequence of ordinal frunctions on $A$ which
is increasing in $<_I$, then if $\l>2^\k$ then $\ov f$ has an exact
upper bound. 

Hint: Define a decreasing sequence of upper bounds of lenth $\k^+$. At
limit $\z<{\k^+}$ prove Lemma \ref{cutstab} using the assumption
$2^\k<\l$ and the pigeon-hole principle. Show that the induction
cannot go on for $\k^+$ steps. Prove lemma
\ref{eub} above using $2^\k<\l$ and the pigeon hole principle. 

See also \cite{jech} for a proof using the Erd\H os-Rado Theorem.
\end{itemize}
\end{exercise}

\section{PCF Theory}

In this Section we begin the study of possible cofinalities of reduced products
of small infinite sets of regular cardinals. We shall see that the set of all
possible cofinalities of such products is well behaved, and is related to a
sequence of ideals over the set.

\begin{definition}
 Let $\lng P,\le\rng$ be a quasi-ordered set, that is $\le$ is a
transitive and reflexive relation over $P$. Write $p<q$ for $p\le
q\wedge q\not \le p$.
 
\begin{enumerate} 
\item A subset $D\su P$ is {\em cofinal} if and only if $\forall p\in
P\exists d\in D\left[p\le d\right]$

\item  The {\em cofinality} $\cf P$ is the least
cardinality of a cofinal subset $D\su P$

\item  A sequence $\ov d= \lng d_i:i<\l\rng$ of elements of $P$ is {\em
increasing cofinal} if and only if $\ov d$ is increasing in $<$ and
$\ran \ov d$ is cofinal

\item $P$ {\em has true cofinality} if and
only if there is an increasing cofinal $\ov d$. The {\em true
cofinality} $\tcf P$ is defined when $P$ has true cofinality and is the
minimal length of an increasing cofinal $\ov d$.
\end{enumerate}
\end{definition}

Every quasi ordered set has cofinality. The disjoint union of
 $\{\om_n:n<\om\}$,
quasi ordered by the union of natural orderings on each $\om_n$, has
cofinality
$\aleph_\om$ and does not have true cofinality.  So the cofinality of $P$
need not be regular, and $P$ may  have no true cofinality.  However,
if $\tcf P$ exists, then $\tcf P= \cf P$ and is regular
(see exercise below).

\begin{definition}
For a set $A$ of regular cardinals let $\pcf A =\{\tcf
\prod A/I :I \text{ is an ideal over} A \text{ and } /prod S/I \text{
has true cofinality}\}$.
\end{definition}

The ordering $\le_I$ was defined for all ordinal functions on $A$, and
it applies in particular for all functions in the product $\prod
A=\{f: f \text { is a function on } A, \forall a\in A\left[f(a)\in
a\right]\}$. Suppose that $\prod A/I$ has true cofinality and fix an
increasing cofinal $\ov f=\lng f_\a:\a<\tcf \prod A\rng$. For every
ideal $I'\supseteq I$ over $A$ the sequence $\ov f$ is increasing
cofinal in $\prod A/{I'}$, and therefore $\tcf\prod A/I=\tcf \prod
A/{I'}$. Since every ideal over $A$ can be extended to a dual of an
ultrafilter, we have proved the following fact:

\begin{fact} 
 $\pcf A=\{\cf \prod A/D: D \text { an ultrafilter over }
A\}$. 
\end{fact}

A product $\prod A/D$ where $U$ is an ultrafilter is linearly
ordered by \L os' Theorem and always has  true cofinality.

\ignore
Why sets and not sequences, and why regular cardinals? For cofinality
purposes it makes no difference replacing a cardinal by its
cofinality; every ideal over a sequence fo cardinal either concentrae
bla bla.
\endignore

We remark that for every set of regular cardinals $A$ we have $A\su
\pcf A$ via principal ultrafilters.

We begin our study of true cofinalities of products. The idea is to
identify which subsets of $A$ are ``too small'' to obtain a regular
$\l$ as true cofinality modulo any ideal. This is the contents of the
very important definition below:

\begin{definition}
Let $A\subseteq  \Reg$ and assume that $\left| A\right| <\min A$.  
For a regular $\lambda $ define $J_{<\lambda }\left[ A\right] =\{B\subseteq
A:B\in D\imply\cf \prod A/D<\lambda $ for all 
ultrafilters $D$ over $A\}$.  
\end{definition}

To verify that $J_{<\l}$ is an ideal see exercise  below.

Since $\cf\prod A/D$ is always regular,  we have

\begin{fact} \label{singisempty}
If $\mu$ is singular then $J_{<\mu^+}=J_{<\mu}$.
\end{fact}

The following Theorem provides information about all products modulo
ideals that extend $J_{<\l}$, and has several important corollaries
for pcf.  The proof is a good warm-up for the next section, whre the pcf Theorem
is proved using $I[l]$. The following proof makes an inessential use of
$I[\l]$ and Theorem \ref{lubexists}. For a proof of it without either
\ref{lubexists} or
$I[\l]$ see [MB] or Lemma 1.5 in chapter 1 of the book.

\begin{theorem}\label{ldirected}
If $\min A>\left| A\right|$ then $\prod A/{J_{<\lambda }}$ is
$\lambda $ -directed.
\end{theorem}

\begin{proof} By induction on $\mu <\lambda $ we show that any
set $F\subseteq \prod A$ with $\left| F\right| =\mu $ has an upper
bound in $\prod A/{J_{<\lambda }}$. 

Let $F\subseteq \prod A$ be given and suppose $\left| F\right| =\mu $
be given.

\medskip
\noindent
\textbf{Case 1}: $\mu $ singular.

Using the induction hypothesis, assume, eithout loss of generality, that
$F=\left\langle f_\alpha :\alpha <\mu \right\rangle$ is increasing in
$<_I$.   Pick a cofinal sequence $\left\langle f_{\alpha _\zeta
}:\zeta <cf\mu
\right\rangle$.  

As $cf\mu <\mu $, the induction hypothesis gives a bound $g\in \prod
A/{J<\lambda }$ to $\left\langle f_{\alpha _\zeta }:\zeta <cf\mu
\right\rangle$ which is also $a$ bound to $\ov f$.

\medskip
\noindent
\textbf{Case 2}:  $\mu$ is regular and $|\mu\cap A|<\aleph_0$.  Let
$g(a)=\sup\{f(a):f\in F\}$ for all $a\ge \mu$ in $A$ and $g(a)=0$ for
$a\in A\cap \mu$. Since $A\cap \mu$ finite, $A\cap \mu\in J_{<\mu}$
and $g\in \prod A$ is an upper bound of $F$.

\medskip
\noindent
\textbf{Case 3}: $\mu$ is regular and $A\cap \mu$ is infinite.
Let $\k:=|A|^+$. Since $|A|<\min A$, we have $\k^+<\mu$.
Use Theorem \ref{beauty} to find $S\subseteq S_{\kappa ^{+}}^\mu$
stationary,
$S\in I\left[ \mu \right]$.   Fix $\ov C=\left\langle c_\alpha
:\alpha <\mu \wedge \cf\a\le\k^{+}\right\rangle$ that demonstrates that
$S\in I\left[\mu\right]$.  Enumerate $F$ as
$\left\langle f_\alpha :\alpha <\mu
\right\rangle$.   Using the induction hypothesis define by induction on
$\a<\mu$ a sequence
$\left\langle h_\alpha :\alpha <\mu \right\rangle$ satisfying

\begin{enumerate}
\item $h_\alpha >_If_\b $ and $h_\alpha >_Ih_\beta $ for all $\beta
<\alpha $  
\item  $\ov h=\left\langle h_\alpha :\alpha <\mu \right\rangle $ obeys 
$\ov C$
\end{enumerate}

This involves taking, at stage $\a$,  an upper bound $g_\a\in \prod A$ of the
set
$\left\{h_\beta :\beta <\alpha \right\} \cup \left\{ f_\beta :\beta <\alpha
\right\}$ and then letting 
$h_\alpha \left( a\right) =\max \left( \left\{ g_\a\left( a\right) \right\}
\cup
\left\{ h_\beta \left( a\right) :\beta \in c_\alpha \right\} \right)$
for all $a\in A\sm \k^{++}$ and letting $h_\a(a)=\a$ for $a\in A\cap
\k^{++}$. 

Use Theorem \ref{lubexists} theorem to find an eub $g$ in $\le_{J_{<\l}}$ of
$\ov h$ (as ordinal function). Because $h_\a(a)<a$ for all $\in A\sm
\k^{++}\in J_{<\l}$,  we may assume, without loss of generality, that
$g\left( a\right) \leq a$ for all $a\in A$.  Let $B=\left\{ a\in
A:g\left( a\right) =a\right\}$.

\begin{claim}
$B\in J_{<\l}$
\end{claim}

\begin{proof}[Proof of Claim]
If $B\notin J_{<\l}$ then we can find and ultrafilter $D\ni B$ over
$A$ so that $D\cap J_{<\lambda }=\emptyset$.  The dual of $D$ extends
$J_{<\l}$ and therefore $\ov h:=\lng h_\a:\a<\mu\rng$ is increasing
also in $<_D$. For every $f\in F$ let $f'$ be obtained from $f$ by
replacing $f\rest( A\sm B)$ by the constant function $0$.  As $B\in
D$, $f=_D f'$. Also $f'<_{J_{<\l}} g$ and therefore, because $g$ is an
eub, there is some $\a<\mu$ so that $f'<_{J_{<\l}} h_\a$. As $D$
extends $J_{<\l}^*$, $f=_Df'<_D h_\a$.  Thus $\ov h$ is increasing
cofinal in $\prod A/D$ demonstrating that $\cf\prod A/D=\mu<\l$. By
definition of $J_{<\l}$ we conclude that $B\in J_{<\l}$.
\end{proof}

Now that we know $B\in J_{<\l}$ redefine $g$ as $0$ on $B$ and it remains
an eub in $\le_{J_{<\l}}$ of $\ov h$, only now belongs to $\prod A$. By the
definition of
$\ov h$, $g$ is an upper bound of $F$.
\end{proof}

We list now several corollaries of the $\l$-directedness Theorem
above. 

\begin{corollary} \label{cofDminl}
If $D$ is an ultrafilter over $A$ then $\cf\prod A/D=\min\{\l:D\cap
J_{<\l^+}\neq \emptyset\}$.
\end{corollary}

\begin{proof}
Suppose $D$ is an ultrafilter over $A$ and let $\l$ be the least such
that $D\cap {J_{<\l^+}}$ is not empty. By the definition of $J_{<\l^+}$ we
have $\cf\prod A/D\le \l$.

Conversely, $D$ extends $J_{<\l}^*$ because $D\cap
 J_{<\l^+}=\emptyset$ and therefore $\le_D$ extends
 $\le_{J_{<\l}}$. Since by Theorem \ref{ldirected} the product
 $\prod A/J_{<\l}$ is $\l^+$-directed, also the product $\prod A/D$
 is $\l^+$-directed, and must therefore have cofinality at least
 $\l$.
\end{proof}

\begin{corollary}
If $\mu$ is a limit cardinal then $J_{<\mu}=\bigcup_{\l<\mu}J_{<\l}$
\end{corollary}\label{minl}

\begin{proof}
Let $J:=\bigcup_{\l<\mu}J_{<\l}$ for some limit cardinal $\mu$. The
inclusion $J\su J_{<\mu}$ is immediate by the definition of
$J_{<\mu}$. 

Conversely, suppose to the contrary that $B\notin J$ and find an
ultrafilter $D\ni B$ over $A$ so that $D\cap J=\emptyset$. Since $D$
extends $J_{<\l}$ for every $\l<\mu$, by Theorem \ref{ldirected} we have
$\cf\prod A/D\ge \l$. Thus $\cf\prod A/D\ge \mu$ and $B\notin
J_{<\mu}$. (Observe that if $\mu$ is  regular limit the equality
$\cf\prod A/D=\mu$ is not out of the question.)
\end{proof}

\begin{corollary}
For every $B\su A$ there is a unique regular $\l$ such that $B\in
J_{<\l^+}\sm J_{<\l}$.
\end{corollary}

\begin{proof}
Since $\cf\prod A/D<|\prod A|^+$ for all ultrafilters $D$ over $A$,
for every $B\su A$ there is some $\mu$, and therefore a minimal $\mu$,
for which $B\in J_{<\mu}$. By the previous Corollary such $\mu$ is not
limit for any $B\su A$; so $\mu=\l^+$ is successor and $B\in
J_{<\l^+}\sm J_{<\l}$. Since $J_{<\mu}$ is increasing in $\mu$, $\l$
is unique. The cardinal $\l$ must be regular by Fact
\ref{singisempty}.
\end{proof}

\begin{corollary}
$\pcf A$ has a last element.
\end{corollary}

\begin{proof}
The previous corollary provides a unique regular $\l$ for which $A\in
J_{<\l^+}\sm J_{<\l}$. Thus $\cf\prod A/D\le \l$ for all ultrafilters $D$
over $A$, and therefore $\l\ge\sup \pcf A$.

Conversely,  let $D$ be an ultrafilter over $A$ so that $D\cap
J_{<\l}\cap D=\emptyset$. Such $D$ exists because $A\notin J_{<\l}$.
The cofinality $\cf\prod A/D$ is at least $\l$ because $D$ extends
$J_{<\l}^*$, and by the previous paragraphs is exactly $\l$, which
means $\l\in \pcf A$. 

Thus $\l=\max \pcf A$.
\end{proof}

\begin{corollary}
If $|A|<\min A$ then $|\pcf A|\le 2^{|A|}$.
\end{corollary}

\begin{proof} By \ref{minl} there is a 1-1 correspondence between $\pcf
A$ and the members of the increasing sequence of ideals $\lng
J_{<\l}:\l\in \pcf A\}$. The length of any increasing sequence of
ideals over $A$ is at most $2^{|A|}$.
\end{proof}

Let us observe the following about the operation $\pcf$ on sets of
regular cardinals. We know that $A\su pcf A$. It is seen (exercise
belos) that $\pcf A\cup B =\pcf A\cup \pcf B$. If we knew that
$\pcf\pcf A=\pcf A$ that this operation would satisfy the axioms
of a topological {\em closure} operation. So let us ask:

\begin{question}
Is $\pcf\pcf A=\pcf A$ for a set of regular cardinals $A$ with $\min
A>|A|$?
\end{question}

This question is answered by a ``yes'' in a universe of set theory
without inaccessible cardinals, as we shall see shortly. But the
important issue is the following: if there is no inaccessible in $\acc
pcf A$ then $\pcf\pcf A=\pcf A$.

\begin{theorem}\label{averaging}
 Suppose $B\su \pcf A$ and $|B|<\min B$. Then $\pcf B\su 
\pcf A$.

\end{theorem}

\begin{corollary} Suppose there is no inaccessible cardinal in $\acc\pcf A$.
Then $\pcf \pcf A=\pcf A$
\end{corollary}

\begin{proof}[Proof of Corollary]

Let  $\l_0:=\sup \pcf A$. Either $\l_0$ is not inaccessible or
it is not an accumulation point of $\pcf A$. In either case, 
$\l_1:=\cf|\pcf A|<\l_0$. Let $B_0=\pcf A\cap(\l_1,\l_0]$. Since $\min
B_0>|B_0|$ we may  use Theorem \ref{averaging} to conclude that $\pcf B_1\su
\pcf A$.  

Continue by induction to define $\l_{n+1}=\cf |\pcf A\cap (\l_n+1)|$ and
$B_{n}:=\pcf A\cap (\l_{n+1},\l_n]$. Because the sequence $\lng
\l_n:n<\om\rng$ is strictly decreasing, it terminates after finitely many
steps (namely, $B_{n+1}=\emptyset$ and $\l_{n+1}=0$ for some $n$) and
thus produces a partition of
$\pcf A$ to finitely many parts $\{B_l:l\le n\}$ so that $\pcf B_l\su \pcf
A$ for all $l\le n$. Since every ultrafilter on $\pcf A$ concentrates on one
of the $B_l$-s, we conclude that $\pcf \pcf A\su \pcf A$.
\end{proof}

In partucular, if there are no inaccessibles in the universe then the
operation of $\pcf$ satisfies the axioms of a topological closure operation.

\begin{proof}[Proof of the Theorem]

 Suppose $B\su \pcf A$, $\min B>|B|$ and fix an ideal $I_\l$  
over $A$ so that $\tcf \prod A/I_\l=\l$ for every $\l\in B$. Suppose that
$\l(*)\in \pcf B$ and that $I(*)$ is an ideal over $B$ so that $\tcf\prod
B/I(*)=\l(*)$. 

Let $I=\left\{X\su A:\{\l\in B:X\notin I_\l\}\in I(*)\right\}$.
We show that $I$ is an ideal over $A$ and
$\tcf\prod A/I=\l(*)$.

The verification that $I$ above is an ideal is immediate. 

For every $\l\in B$ fix an increasing cofinal sequence $\ov f^\l=\lng
f^\l_\a:\a<\l\rng$ in $\prod A/I_\l$. 

Define for every function $g\in \prod B$ a function $G(g)\in \prod A$ as
follows: since $\prod A/J_{<|B|^+}$ is $|B|$-directed by Theorem
\ref{ldirected}, there is a bound $f\in \prod A/J_{<|B|^+}$ to the set
$\left\{f^\l_{g(\l)}:\l\in B\right\}$.  Let
$G(g)$  be such a bound.

Fix an increasing cofinal sequence $\left\lng
g_\a:\a<\l\left(*\right)\right\rng$ in
$\prod B/I(*)$. We argue that for every $f\in \prod A$ there is some
$\b<\l(*)$ such that $f<_IG(g_\a)$ for all $\a\in (\b,\l(*))$. Suppose that
$f\in
\prod A$ is given. Let $F(f)\in \prod B$ be defined by
$F(f)(\l)=\min\{\a<\l:f<_{I_\l} f^\l_\a\}$ for $\l\in B$. This is well
defined because $\lng f^\l_\a:\a<\l\rng$ is increasing cofinal in $\prod
A/I_\l$ for all $\l\in B$. Let $\b<\l(*)$ be the least for which
$F(f)<_{I(*)} g_\b$. If $\a\in (\b,\l(*))$ then $F(f)<_{I(*)}g_\a$
and thus $f<_{I_\l}f^\l_{g_\a(\l)}$ for all but a set in $I(*)$ of
$\l\in B$. This means that $f<_I G(g_\a)$ by the definition of $I$.

The above implies that $\lng G(g_\a):\a<\l(*)\rng$ is cofinal in $\prod A/I$
and that for some club $E\su \l(*)$ the sequence $\lng G(g_\a):\a\in E\rng$
is increasing in $<_I$. That establishes $\tcf\prod A/I=\l(*)$. 

Thus $\pcf B\su \pcf A$ and the proof is done.
\end{proof}

\section{The existence of generators for pcf}\label{gens}

In this section we will prove that for every $\l\in pcf A$, $A\su \Reg$
with $\min A>|A|$, there is a set $B_\l\su A$ that {\em generates} the
ideal $J_{<\l^+}$ from $J_{<\l}$, that is $J_{<\l^+}=J_{<\l} + B_\l$,
the minimal ideal extending $J_{<\l}$ which contains $B_\l$. Such
$B_\l$ are called {\em generators} for $\pcf A$. The generator $B_\l$
is unique only up to ${J_{<\l}}$, of course. The theorem asserting the
existence of generators for pcf (Theorem \ref{genexists}) is called
{\em the pcf Theorem}.

\witb

The book defines in Chapter 1 a $\l\in \pcf A$ for which a generator
exists as {\em normal}. In this terminology, we are proving that every
$\l$ is normal. In \cite{Sh506} the pcf Theorem is re-proved under
more general assumptions.

The following Lemma has particular importance:

\begin{lemma}\label{tcfjl}
If $\l\in \pcf A$ and $B\in J_{<\l^+}\sm J_{<\l}$ then $\tcf \prod
B/J_{<\l}=\l$.
\end{lemma}

\begin{proof}
First let us observe that if $\l\cap A$ is finite there is little to
prove, as $J_{<\l^+}=A\cap \l^+$ and $J_{<\l}=A\cap \l$, and thus
$J_{<\l^+}\sm J_{<\l}=\{\l\}$. So we assume
that $\l>|A|^{++}$ and fix a stationary set $S\su S^\l_{|A|^+}$ in
$I[\l]$ together with a sequence $\ov C$ witnessing this.

Define $J:=J_{<\l}\cup \{B\in J_{<\l^+}\sm J_{<\l}:\tcf
B_{J_{<\l}}=\l\}$.

\begin{claim}
$J$ is an ideal. 
\end{claim}

\begin{proof}[Proof of Claim]

 Suppose that $B_1,B_2\in J_{<\l^+}\sm J_{<\l}$ are in $J$ and fix
increasing cofinal sequences $\ov f^1,\ov f^2$ in $\prod B_1/J_{<\l}$,
$\prod B_2/J_{<\l}$ respectively to demonstrate this. Let $f_\a=\max
\{f^1_\a,f^2_\a\}$. The sequence $\ov f$ is increasing in
$\le_{J_{<\l}}$ and cofinal on $B_1\cup B_2$, and this is enough for
$B_1\cup B_2\in J$. The cases in which one of $B_1,B_2$ belongs to
$J_{<\l}$ are immediate. 
\end{proof}

If $J=J_{<\l^+}$ then we are done. Suppose that $B\in J_{<\l^+}\sm J$ and
we shall derive a contradiction. 

Fix and ultrafilter $D\ni B$ over $A$ so that $J\cap D=\emptyset$. 
By Corollary \ref{cofDminl} we conclude that $\cf\prod A/D=\l$. Fix,
then, an increasing cofinal sequence $\lng f_\a:\a<\l\rng$ in $\prod
A/D$. By induction on $\a<\l$ choose $h_\a\in \prod B$ so that $\lng
h_\a:\a\in\l\rng$ is increasing in $<_{J_{<\l}}$ and
$h_\a>_{J_{<\l}}f_\b$ for all $\b<\a$. The definition is possible by
$\l$-directedness modulo $J_{<\l}$ (Corollary \ref{ldirected}).
Furthermore, make $\ov h$ obey $\ov C$.

Let $g$ be an $\eub$ of $\ov h$ so that $g(a)\le a$ for all $a\in A$
and denote $B_g:=\{ a\in A:g(a)=a\}$. Since $\ov h$ is increasing
cofinal below $g$ by definition of an $\eub$, we have that $\tcf \prod
B_g/J_{<\l}=\l$, and therefore $B_g\in J$. On the other hand, if $B_g$
were not in $J$, then $A\sm B_g\in D$.  Since $g\rest (A\sm B_g)\cup
0\rest B_g\in \prod A$, we obtain that $\ov h$ and consequently $\ov
f$ are bounded in $\prod A/D$ --- contradiction.  So $B_g\in D$. This
contradicts $D\cap J=\emptyset$.
\end{proof}

A generator is a maximal set with respect to $\su_{J_{<\l}}$. To find
it we will construct, in the proof of the pcf theorem below, an
increasing sequence in $\subsetneqq_{J_{<\l}}$.  To be able to use of
$I[\l]$ for the purpose of mixing $\subsetneqq_{J_{<\l}}$ with $\su$,
we will use a special kind of increasing cofinal obedient sequences,
which we call {\em humbly obedient}.

\begin{definition} Let $\l\in \pcf A\sm |A|^{+++}$, and let $\ov C$
witness that a stationary $S\su S^\l_{\k^+}$ is in $I[\l]$. A sequence
$\ov f=\lng f_\a:\a\in \l\rng$ {\em humbly obeys} $\ov C$ if and only
if $\ov f$ is increasing in $<_{J_{<\l}}$, $\ov f$ obeys $\ov C$ and
$f_\d=\sup\{f_\b: \b\in c_\d\}$ for every $\d\in S$.
\end{definition}

\noindent \textbf{Explanation}: To obey $\ov C$ means that $f_\a\ge
\sup \{f_\b:\b\in c_\a\}$. To humbly obey means that
for $\d\in S$, $f_\d=\sup\{f_\b:\b\in c_\d\}$, namely $f_\d$ is the
{\em minimal} function allowed by obedience when $\d\in S$. Notice
that while to obey was defined for all sequences we reserve humble
obedient for increasing sequences in $\prod A/{J_{<\l}}$ of length
$\l$.

Humbly obedient sequences form a proper subclass of obedient
sequences. Let us show that it is not an empty subclass by showing we
can find for every increasing sequence $\ov f$ a humbly obedient
$\ov h$ which bounds at least what $\ov f$ bounds.

Fix now $\l\in \pcf A\sm |A|^{+++}$, $S$ and $\ov C$ as in the
definition for the rest of this section.

\begin{fact} \label{humbdense} For every sequence $\ov f=\lng
f_\a:\a<\l\rng$ which is increasing in $<_{J_{<\l}}$ there exists a
humbly obedient sequence $\ov h=\lng h_\a:\a<\l\rng$ such that $\ov h$
is increasing in $<_{J_{<\l}}$ and $f_\a\le h_\a$ for all $\a\in \l\sm
S$.  \end{fact}

\begin{proof} By induction on $\a<\l$ define $h_\a$. For $\a=\d\in S$
define $h_\a=\sup\{h_\b:\b\in c_\d\}$. For $\a\in \l\sm S$ find, using
$\l$-direc tedness, a bound $g_\a\in \prod A$ of
$\{f_\b:\b<\a\}\cup\{h_\b:\b<\a\}$ and define
$h_\a(a)=\sup\left(\{g_\a(a)\}\cup \{f_\a(a)\} \cup
\{h_\b(a)+1:\b\in c_\a\}\right)$.
\end{proof}

Let us define now a partial ordering on the class of humbly obedient
sequences.

\begin{definition}
If $\ov f^1,\ov f^2$ are humbly obedient, write $\ov f^1\le \ov f^2$
if and only if $f^1_\a\le f^2_\a$ for all $\a<\l$.
\end{definition}

\begin{fact}
If $\ov f^1,\ov f^2$ are humbly obedient then $\ov f^1\le \ov f^2$ if
and only if $f^1_\a\le f^2_\a$ for all $\a\in \l\sm S$.
\end{fact}

\begin{proof} One direction follows because $\l\sm S\su \l$. For the
other direction suppose that $f_\a^1\le f^2_\a$ for all $\a\in \l\sm
S$ and fix a least $\d\in S$ for which $f^1_\d\not\le f^2_\d$. Since
for every $\b\in c_\d$ we have that $f^1_\b\le f^2_\b$ and
$f^2_\d=\sup\{f_\b:\b\in c_\d\}$, $f^1_\d=\sup\{f^1_\b:\b\in c_\d\}$,
contradiction follows.
\end{proof}

Since every humbly obedient sequence $\ov f$ is in particular an
increasing in $<_{J_{<\l}}$ obedient sequence of ordinal functions, it
has, by
Theorem \ref{lubexists}, an exact upper bound. Choose eub
 $g_{\ov f}$ for every humbly obedient $\ov f$, such, without loss of
generality, that $g_{\ov f}(a)\le a$ for all $a\in A$. Let us also denote $B_{\ov
f}:=\{a\in A: g_{\ov f}(a)=a\}$. Since $\ov f$ is cofinal below $g$, we obtain that
$\tcf B_{\ov f}/{J_{<\l}} =\l$ and therefore $B_{\ov f}\in J_{<\l^+}\sm
{J_{<\l}}$ for every humbly obedient $\ov f$. 

We show next that every
set $B\in J_{<\l^+}\sm {J_{<\l}}$ is equivalent mod ${J_{<\l}}$ to a
$B_{\ov f}$ for some humbly eobedient $\ov f$, and that if $\ov f^1\le
\ov f^2$ then $B_{\ov f^1}\su_{J_{<\l}} B_{\ov f_2}$.

\begin{lemma}  \label{rep}
\begin{enumerate}
\item Suppose $B\in J_{<\l^+}\sm {J_{<\l}}$. Then there is a
humbly obedient $\ov f$ such that $B=_{J_{<\l}} B_{\ov f}$.
\item If $\ov f^1\le \ov f^2$ are humbly obedient sequences, then
$B_{\ov f_1}\su_{J_{<\l}} B_{\ov f^2}$
\end{enumerate}
\end{lemma}

\begin{proof}
Extend the ideal ${J_{<\l}}$ by the set $A\sm B$. Thus we may ignore
any coordinalte outside of $B$ and work in either $\prod A/{J_{<\l}}$
or $\prod B/{J_{<\l}}$.

Using lemma \ref{tcfjl} we fix an increasing cofinal seqeunce of length
$\l$ in $\prod B/{J_{<\l}}$ and replace it with a humbly obedient
increasing cofinal $\ov f$, which exists by lemma
\ref{humbdense}. Let $g$ be an eub of $\ov f$ and $B_{\ov
f}=B_g=\{a\in A:g(a)=a\}$. Since we assume that $\ov f$ is constantly
zero outside of $B$, we may assume that $g$ is zero outside of $B$ and
therefore $B_g\su B$.

Let now $C=B\sm B_g=\{a\in B:g(a)\in a\}$. The sequence
$\lng f_\a\rest C:\a<\l\rng$ is increasing cofinal in $\prod
C/{J_{<\l}}$ and therefore $g\rest C<_{J_{<\l}} f_\a$ for some
$\a<\l$. Since $f_\a<_{J_{<\l}} g$ we conclude that $C\in {J_{<\l}}$.

Thus $B=_{J_{<\l}} B_g$.

To prove 2, observe that since $f^1_\a\le f^2_\a$ for all $\a<\l$ and
$\ov f^1$ is cofinal $\prod B_{\ov f^1}/{J_{<\l}}$, also $\ov f^2\rest
B_{\ov f^1}$ is cofinal in $\prod B_{\ov f^1}/{J_{<\l}}$, and
therefore $B_{\ov f^1}\su_{J_{<\l}} B_{\ov f^2}$.
\end{proof}

If a generator $B_\l$ for $J_{<\l^+}$ over ${J_{<\l}}$ does exist,
then $J_{<\l^+}/{J_{<\l}} $ is a boolean algebra, namely has a maximal
element $B_\l$ and is in particular $\k$-directed for all $\k$. The
following lemma, which is needed to prove the pcf theorem, asserts a
particular case of directedness: that $J_{<\l^+}\sm {J_{<\l}}$ is
$|A|^+$-directed. It  says actually slightly more:  that the
relation $\le$ on humbly obedient sequences is $|A|^+$-directed.

\begin{lemma}  \label{humbdirected}
Suppose that   $\z(*)<|A|^+$. 
\begin{enumerate} 
\item If  $\{ \ov f^\z:z<\z(*)\}$ is a set of
humbly obedient sequences, then there
exists a humbly obedient sequence $\ov f$ so that $f_\z\le f$ for all
$\z<\z(*)$.
\item If $\{ B_\z:\z<\z(*)\}\su J_{<\l^+}$  then there is a set $B\in
J_{<\l^+}$ so that $B_\z\su_{J_{<\l}} B$.
\end{enumerate}
\end{lemma}

\begin{proof}
Because every $B\in J_{<\l^+}$ is equivalent to $B_{\ov f}$ for some
humbly obedient $\ov f$, and $\ov f_\z\le \ov f$ implies $B_{\ov
f^\z}\su_{J_{<\l}} B_{\ov f}$, and $B_{\ov f}\in J_{<\l^+}\sm
{J_{<\l}}$ for all humbly obedient $\ov f$, it is enough to prove the
first item in the lemma.

Suppose that $\{\ov f^\z:\z<\z(*)\}$ are humbly obedient and
$\z(*)<|A|^+$. Define by induction on $\a<\l$ an increasing in
$<_{J_{<\l}}$, humbly obedient sequence $\lng f_\a:\a<\l\rng$ such
that $f^\z_\a\le_{J_{<\l}} f_\a$ for all $\z<\z(*)$, $\a<\l$. At the
induction step you need to immitate the proof of
\ref{humbdense} above, with the only change of taking an upper bound
$g_\a$ over a union of $\z(*)\times \a$ many functions.
\end{proof}

The following lemma is the key to the proof of the pcf theorem:

\begin{lemma}\label{aggreeonclub}
If $\ov f^1\leq \ov f^2$ are humbly obedient sequences and $B:=B_{\ov f^1}\in
J_{<\lambda^+ }\sm J_{<\l}$ then there is a club $E\subseteq \lambda $ such
that  
\begin{equation}
\delta \in S\cap E\Rightarrow f_\alpha ^1\rest B
=_{J_{<\lambda }}f_\delta ^2\rest B
\end{equation}
\end{lemma}

\begin{proof}  The sequence $\ov f^1$ is increasing
cofinal in $J_{<\lambda }$ on $B$ by the definition of $B$.   Thus, for every
$\alpha <\lambda$ there is some $\beta \left( \alpha \right) <\lambda $ such
that $f_\alpha^2<_{J_{<\lambda }}f_\alpha ^1$. Let $E\subseteq \lambda$ be a
club closed under $\alpha \mapsto \beta \left( \alpha \right) $ (that is,
$\delta
\in E\,\wedge \alpha <\delta \Rightarrow \b(\alpha) <\delta$).  
Suppose that $\delta \in E\cap S$. We show that $f_\delta ^1\rest
B=_{J_{<\lambda }}f_\delta ^2\rest B$.   We already know that $f_\delta ^2\geq
f_\delta ^1$.   

Conversely, let 
$C=\left\{ a\in B:f_\delta ^1\left( a\right) <f_\delta ^2\left( a\right)
\right\}$.

For every $a\in A$ we have 
$f_\delta ^2( a) =\sup \left\{f_\alpha ^2(a) :\alpha
\in c_\delta \right\}$, and therefore
for every $a\in C$  we can find an index $\alpha _a\in c_\delta $ such that 
$f^2_{\alpha _a}(a) >f_\delta ^1\left( a\right)$.   Such an index
must exist because $f_\delta ^2( a)=\sup\{f_\a(a):a\in c_\d\}
>f_\delta ^1\left( a\right)$.   As $\otp c_\d=|A|^+$, while $|C|\le |A|$ we
can find  
$\gamma \in\nacc c_\delta $ which is greater than $\a_a$ for all $a\in C$.  By
coherence at
$\gamma$ it follows that
$f^2_\gamma(a)>f^2_\a(a)>f^1_\d(a)$ for all $a\in C$. 
But because $\d\in E$ there is some $\b<\d$ for which
$f^2_{\gamma}<_{J_{<\l}}f^1_\b$. The  set $C':=\{a\in
A:f^2_\gamma(a)<f^1_\b(a)<f^1_\d(a)\}$ cannot  meet  $C$, since
$a\in C'\cap c\imply f^2_\gamma(a)<f^1_\b(a)<f^1_\d(a)<f^2_\gamma(a)$, which is
absurd. As $C'\cup A\sm B$ is measure 1, it follows that $C\in J_{<\l}$ and
$f^2_\d\rest=_{J<\l}f^1_\d\rest B$
v\end{proof}

We have all the facts  about humbly obedient needed to prove the pcf Theorem.
Let us devote a few words to the description of what is going on in
this proof, compared to the proof of the lub theorem (Theorem
\ref{lubexists} above). 

In the proof of the lub theorem we were able to mix the relations
$\su$, which appears locally, and $\subsetneqq_I$, which occurs on a
club, and thus obtain a contradiction from a sequence of lenght
$|A|^+$ of subsets of $A$. That proof had a ``one-dimensional''
geometry.

This proof uses a similar idea, but is  ``two-dimensional'', by which
we mean that the increasing sequence of subsets of $A$ of length
$|A|^+$ is indexed by pairs of ordinals.   The
phenomenon that reflects on a club is not increasing
in $\su_{J_{<\l}}$, but equality modulo ${J_{<\l}}$

\begin{theorem}\label{genexists}
Suppose $A\subseteq \Reg$ is infinite and $\left| A\right|  <\min A$.
Then for every $\lambda \in \pcf A$ there is a set $B_\lambda
\subseteq A$ such that
\[
J_{<\lambda ^{+}} =J_{<\lambda } +B_\lambda
\]
\end{theorem}

\begin{proof} Since for every $\l\in \pcf A$ for which $\l\cap A$ is
finite the theorem does not say more than $\l^+\cap A=\l \cap A\cup
\{\l\}$, the Theorem is true for those $\l$. 

Suppose, then, that $\lambda \in \pcf A\sm |A|^{++}$ and fix a
stationary $S\subseteq S^{\l}_{\k^+}$ in $I\left[ \lambda
\right]$ and a sequence $\ov C$ witnessing  this. We shall find
a generator for $J_{<\l^+}$ over ${J_{<\l}}$.

 By induction on $\zeta <|A|^{+}$ we define humbly obedient sequences $\ov
f^\z$ with $B_\z:=B_{\ov f^\z}$ so that
\begin{itemize}
\item
$\xi<\z<|A|^+\imply \ov f^\xi\le \ov f^\z$
\item $\xi<\z<|A|^+\imply B_\xi\subsetneqq_{J_{<\l}}  B_\z$
\end{itemize}

Suppose that $\z<|A|^+$ and that $\ov f^\xi$ and $B_\xi:=B_{\ov
f^\xi}$ are defined for all $\xi<\zeta$ and satisfy the conditions
above. By $|A|^+$-directedness (Lemma
\ref{humbdirected} above) there is a set $B\in J_{<\l^+}$ such that
$B_\xi\su_{J_{<\l}} B$ for all $\xi<\z$. If $B$ generates $J_{<\l^+}$,
we are done. Else, there is a set $B'$ so that $B\subsetneqq_{J_{<\l}}
B' \in J_{<\l^+}$. Use Lemma \ref{rep} to fix a humbly obedient
sequence $\ov f$ so that $B'=B_{\ov f}$. Now let $\ov f_\z$ be
provided by Lemma \ref{humbdirected} so that $\ov f\le \ov f^\z$ and
$\ov f_\xi\le
\ov f^\z$ for all $\xi<\z$. Let $B_\z=B_{\ov f^\z}$ and it follows that
$B_\xi\subsetneqq_{J_{<\l}} B_\z$ for all $\xi<\z$.

Suppose to the contrary that $\ov f^\z$ and $B_\z$ are defined for all
$\z<|A|^+$ and satisfy that $\ov f^\z$ is increasing in $\le$ and
$B_\z$ increasing in $\subsetneqq_{J_{<\l}} $. 

For every pair $\xi<\z<|A|^+$ there is, by Lemma \ref{aggreeonclub}, a
club $E_{\xi,\z}\su \l$ such that $\d\in S\cap E_{\xi,\z}\imply
f^\xi_\d\rest B_\xi=_{J_{<\l}} f^\z_\d\rest B_\xi$. Since $\ov f^\z$
is increasing cofinal on $B_\z\sm B_\xi$ and $g_\xi\rest B_\z\sm
B_\xi\in \prod B_\z\sm B_\xi$ by definition of $B_\xi$, there is some
$\a<\l$ for which $g_\xi\rest (B_\z\sm B_\xi)<_{J_{<\l}} f^\z_\a$. Thus,
by subtracting $\a$ from $E_{\xi,\z}$ we may assume that $f^\xi_\d\rest
B_\z\sm B_\xi<_{J_{<\l}} f^\z_\d$ for all $\d<\l$ and $\z\in E$.

Let $E=\bigcap_{\xi<\z<|A|^+}E_{\xi,\z}$.  Since $|A|^+<\l$, the
intersection $E$ is a club of $\l$. Fix $\d\in S\cap E$. For all
$\xi<\z<|A|^+$ we have
\begin{equation}\label{splitting} 
f^\z_\d\rest B_\xi=_{J_{<\l}} f^\xi_\d\rest B_\xi \text{  and }  f^\xi_\d\rest
(B_\z\sm B_\xi) <_{J_{<\l}} f^\z_\d\rest (B_\z\sm B_\xi)
\end{equation}

Picture?????

We do the last bit of the proof twice, for no real reason other than
to emphasize the similarity to the proof of the lub theorem. 

\noindent
\textbf{First run}:
 For every $a\in A$,the sequence $\lng f^\z_\d(a):\zeta<|A|^+\rng$ is
increasing in $\le$ because $\lng \ov f_\z:|z<|A|^+\rng$ is increasing
in $\le$. Therefore, for every $a\in A$ there is a club $C_a\su |A|^+$ on which
$\lng f^\z_\d(a):\zeta<|A|^+\rng$ is either constant or is strictly
increasing (that is, increasing in $<$).

Let $C=\bigcap_{a\in A}C_a$.
 Since $|A|<|A|^+$, this intersection is a
club of $|A|^+$.

Choose  $\xi<<\rho$ in $C$ and $\z=\zeta(\xi)<\rho$. Since $\d\in E\su
E_{\zeta,\rho}$ we know that $f^\z_\d(a)=f^\rho_\d(a)$ for all $a\in
B_\z$ but a measure zero set. On the other hand, the set $\{\a\in
B_\z\sm B_\xi:f^\xi(a)<f^\z(a)\}$ is positive. That allows us to find
$a\in B_\z\sm B_\xi$ on which both relations occur: namely
$f^\xi_\d(a)<f^\z_\d(a)=f^\rho_\d(a)$. This is impossible, though,
because $\xi,\z,\rho\in C\su C_a$ and therefore $\xi<\z<\rho$ implies
that either $f^\xi_\d(a)=f^\z_\d(a)=f^\rho_\d(a)$ or
$f^\xi_\d(a)<f^\z_\d(a)<f^\rho_\d(a)$.

\noindent
\textbf{Second run}: For every $\xi<\k$ define $X_{\z,\xi}:\left\{a\in
A:f^\z_\d(a)=f^\xi_\d(a)\right\}$ for every $\z<|A|^+$. Since $\lng \ov
f_\d^\z:\z<|A|^+\rng$ is increasing in $\le$, the sequence $\lng
X_{\z,\xi}:\z<|A|^+\rng$ is decreasing in $\su$ and stabilizes at
some $\zeta(\xi)$. Let $C\su |A|^+$ be a club of $|A|^+$ which is
closed under $\xi\mapsto
\z(\xi)$. 

If $\xi<\z\le\xi'<\z'$ are in $C$ then $ X_{\z,\xi}=X_{\z',\xi}\su
X_{\z',\xi'}$. The first equality is because $\z,\z'\in C$; the second
inclusion is becaus $f^\xi_\d\le f^{\xi'}_\d\le f^{\z'}_\d$. Thus, for all $
\xi<\z\le\xi'<\z'$ in $C$:

\begin{equation} \label{crucialequation}
X_{\z',\xi'}\su X_{\z,\xi}
\end{equation}

On the other hand, if $\xi<\z\le\xi'<\z'$ are in $C$ then
$B_{\xi'}\su_{J_{<\l}} X_{\xi',\z'}$ by \ref{splitting} and by the
same also $B_\xi'\sm B_\z\cap X_{\xi,\z}=_{J_{<\l}}\emptyset$. Hence:

\begin{equation}\label{reallycrucial}
X_{\z',\xi'}\not=_{J_{<\l}}  X_{\z,\xi}
\end{equation}

Combining \ref{crucialequation} with \ref{reallycrucial} and denoting
by $\xi':=\min \{C\sm (\xi+1)\}$ we get that the sequence $\lng
X_{\xi,\xi'}:\xi\in C\rng$ is a strictly ioncreasing sequence of
subsets of $A$ of length $|A|^+$ - a contradiction. 
\end{proof}

\end{document}